\newif\ifsmfart
\IfFileExists{smfart.cls}
   {\documentclass[12pt,english]{smfart}
    \IfFileExists{smfenum.sty}{\usepackage{smfenum}}{}
    \usepackage{bull}
    \smfarttrue}
{\message{^^J*** It would be better to typeset this file with smfart.cls ***^^J^^J}

\documentclass[12pt]{amsart}}

\usepackage{amssymb}
\usepackage{epsfig}
\usepackage{graphicx}
\numberwithin{equation}{section}

\input xy
\xyoption{all}

\advance\headheight 2pt
\calclayout
\allowdisplaybreaks[3]

\theoremstyle{plain}
\newtheorem{prop}[subsubsection]{Proposition}
\newtheorem{Theo}[section]{Theorem}
\newtheorem{thm}[subsubsection]{Theorem}
\newtheorem{coro}[subsubsection]{Corollary}

\newtheorem{lemm}[subsubsection]{Lemma}
\newtheorem{assu}[subsubsection]{Assumption}
\newtheorem{defn}[subsubsection]{Definition}
\newtheorem{sublemm}[subsubsection]{Sublemma}

\theoremstyle{definition}

\theoremstyle{remark}
\newtheorem{rem}[subsubsection]{Remark}

\newtheorem{exam}[subsubsection]{Example}

\newtheorem{nota}[subsubsection]{Notations}

\newcommand{\cF}{\mathcal{F}}

\newcommand{\tf}{\tilde{f}}

\newcommand{\ovl}{\overline}

\newcommand{\la}{\lambda}
\newcommand{\al}{\alpha}

\newcommand{\Q}{\Bbb Q}
\newcommand{\Z}{\Bbb Z}

\def\cA{{\mathcal A}}

\def\cF{{\mathcal F}}

\def\cO{{\mathcal O}}

\def\I{{\mathrm I}}
\def\cI{{\mathcal I}}

\def\cR{{ R}}
\def\cT{{\mathcal T}}

\def\char{{\rm char}}

\newcommand{\G}{\Gamma}

\newcommand{\Hom}{{\rm Hom}}
\newcommand{\Ker}{{\rm Ker}}

\newcommand{\Gal}{{\rm G}}
\newcommand{\K}{{\rm K}}

\def\ra{\rightarrow}
\def\da{\downarrow}

\def\A{{\mathbb A}}

\def\P{{\mathbb P}}
\def\Q{{\mathbb Q}}

\def\Z{{\mathbb Z}}

\def\N{{\mathbb N}}

\def\GL{{\rm GL}}

\def\rk{{\rm rk}}

\def\AF{{\mathcal A}{\mathcal F}}

\def\LF{{\mathcal L}{\mathcal F}}

\makeatother
\makeatletter

\author{Fedor Bogomolov}
\address{Courant Institute of Mathematical Sciences, N.Y.U. \\
 251 Mercer str. \\
 New York, NY 10012, U.S.A.}
\email{bogomolo@cims.nyu.edu}

\author{Yuri Tschinkel}
\address{Department of Mathematics \\
         Princeton University\\ 
         Fine Hall, Washington Road\\
         Princeton, NJ 08544-1000,  U.S.A.}
\email{ytschink@math.princeton.edu}

\title[Galois groups of function fields]{Commuting 
elements\\
in Galois groups  of function fields}

\begin{document} 
 
\date{\today}



\begin{abstract}
We study the structure of abelian subgroups of
Galois groups of function fields. 
\end{abstract}

\maketitle

\tableofcontents

\setcounter{section}{0}
\section*{Introduction}
\label{sect:introduction}

\subsection*{Setup}

We fix a prime number $p$. Let $k$ be
a field of characteristic $\neq p$ and $2$. 
Assume that $k$ does not admit finite extensions of degree divisible
by $p$. Let $K=k(X)$ be the field of functions of an
algebraic variety $X$ defined over $k$ and 
$\Gal_K$ the Galois group
of a separable closure of $K$.
The principal object of our study 
is the group $\G=\G_K$  - the (maximal) pro-$p$-quotient 
of the kernel  
$\Ker(\Gal_K\ra \Gal_k)$.

\subsection*{Main theorem}

This paper contains a proof of the main
theorem from  \cite{B2} describing (topologically) 
noncyclic subgroups in the abelianization  
$$
\G^a:= \G/[\G,\G]
$$    
which can be lifted to abelian subgroups of 
$$
\G^c:=\G/[[\G,\G],\G].
$$

Let $\nu$ be a valuation of $K$. We denote
by $K_{\nu}$ the completion of 
$K$ with respect to $\nu$, by $\G^a_{\nu}$ the reduced
valuation group and by $\I^a_{\nu}$ the abelian inertia group
of $\nu$ (see Sections~\ref{sect:valuation-notations} and 
\ref{sect:inertia-group} for the definitions).

\begin{Theo}\label{thm:main}
Let $F$ be a noncyclic 
subgroup of $\G^a$. 
Suppose that it can be lifted to an abelian
subgroup of $\G^c$.  
Then there exists a nonarchimedean valuation
$\nu$ of $K$ such that

$\bullet$\,\, $F$ is contained in the abelian
reduced valuation group $\G^a_{\nu}\subset \G^a$ 
(standard valuation group
if the residue field of $K_\nu$ 
has characteristic $\neq p$);

$\bullet$\,\, $F$ contains a subgroup $F'$ such that
$F'\subset \I_v^a$ and $F/F'$ is topologically cyclic. 
\end{Theo}

It is easy to see that abelian groups 
satisfying the conditions of the
theorem can be lifted to 
abelian subgroups of $\Gal_K$
itself. Thus Theorem~\ref{thm:main} shows that 
there are no
obstructions to the lifting of abelian subgroups
of $\Gal_K$ beyond the first nontrivial level.

\

Our paper is a contribution to 
the ``anabelian geometry'' program, initiated
by Grothendieck. For other results in this direction  
we refer to \cite{pop1}, \cite{pop3}, \cite{mochizuki}.

\subsection*{Structure of the proof}
\label{sect:leitfaden}

By Kummer theory, the elements of $\G^{a}$ can be interpreted
as $k^*$-invariant $\Z_p$-valued logarithmic 
functions on $K^*$. 
The quotient space $K^*/k^*$ has a 
natural structure of an infinite dimensional
projective space over $k$, 
denoted by $\P(K)$.   
Consider a pair of 
(nonproportional) elements $f_1,\,f_2$ of $\G^{a}$.
They define a map
\begin{equation}\label{eqn:varphi}
\begin{array}{ccc}
\varphi : &  \P(K)  \ra  & \A^2(\Z_p)\\
          &   (v,v')\mapsto & (f_1(v),f_2(v')).
\end{array}
\end{equation}
If $f_1,f_2$ lift to a pair of commuting 
elements in $\G^c$ then the restrictions of 
the corresponding functions to any projective
line $\P^1\subset \P(K)$ 
are linearly dependent
modulo constant functions (see Proposition~\ref{lemm:cp}).
Thus every projective line in $\P(K)$ maps
into an affine line in $\A^2$.  
This - together with the logarithmic property of $f_1$ and $f_2$ -
imposes very strong conditions on the 2-dimensional 
subspace  they span in the space of functions on $\P(K)$. 
Namely, this subspace contains 
a special nonzero function which we
call an abelian flag function (an AF-function);
it corresponds to an inertia element of some reduced
valuation subgroup (see 
Section~\ref{sect:inertia-group} and 
Section~\ref{sect:fg} for the definitions). 
The main problem is to prove the 
existence of this AF-function. 

\

In Section~\ref{sect:fg} we 
define AF-functions and study them on abelian groups of
ranks 2 and 3. 
Let $f$ be a function  
on a vector space $V$ with values in a set $S$, 
$V'\subset V$ a subspace and $f_{V'}$ the restriction of
$f$ to $V'$. 
A series of reductions  leads to the following 
criterium: $f$ is an (invariant) AF-function on $V$ iff 
for all 3-dimensional subspaces $V'$ the function
$f_{V'}$ is an (invariant) AF-function on $V'$ (\ref{lemm:restr3}).
Reduction to 2-dimensional subspaces is more problematic.
For fields $k$ is of characteristic $\char(k)>2$ the reduction to 
dimension 2 can be established without the use of the logarithmic property 
leading to an easier proof of this case 
of the main theorem. 
A similar statement for fields $k$ of characteristic
zero requires the logarithmic property.
The corresponding proofs are in Section~\ref{sect:checking}.

\

The proof of the main theorem proceeds by 
contradiction. We assume that 
the $\Z_p$-span $\langle f_1,f_2\rangle_{\Z_p}$ does not contain an 
AF-function. 
The reductions and the logarithmic property  
imply that there exists a 3-dimensional $V\subset K$, 
two nonproportional 
functions $f_1',f_2'\in \langle f_1,f_2\rangle_{\Z_p}$ and 
a map $h':\Z_p\ra \Z/2$ such that each of  
the functions $h'\circ f_1', h'\circ f_2', h\circ f_1'+ h\circ f_2'$
fails to be AF. 
In Section~\ref{sect:geometry} we find a 
contradiction to this claim.

\

\

{\bf Acknowledgements.} 
The first author was partially supported by the NSF. 
The second author was partially supported by the NSA.

\section{Classes of functions}
\label{sect:fg}

In this section we define certain classes of functions
on abelian groups and vector spaces which will be used 
in the proofs.  

\subsection{Notations}
\label{not:AS}

Denote by $\Z_{(0)}=\Z\setminus 0$ and by $\Z_{(q)}$ the
set of all integers coprime to $q$. 
Let ${\mathcal A}_0$  (resp. ${\mathcal A}_q$) be the set of  
torsion free abelian groups (resp. vector spaces 
over the finite field ${\bf F}_q$, 
where $q$ is a prime number $\neq 2$). 
We denote by $\rk(A)\in \{1,...,\infty\}$ 
the minimal number of generators of $A$ (as an abelian group). 
An element $a\in A$ (with $A\in {\mathcal A}_0$) 
is called primitive
if there are no $a'\in A$ and $n\in \N$ such that $a= na'$. 
We denote by $\langle a_1,...,a_n\rangle$ 
the subgroup 
generated by $a_1,...,a_n$ and similarly by
$\langle B\rangle$ the subgroup generated by elements of 
a subset  $B\subset A$.
We denote by ${\cF}(A,S)$
the set of functions on $A$ with values in a set $S$.
We will say that a function 
$f\in \cF(A,S)$ is induced from $A/nA$
if for all primitive $x$ and all 
$y\in A$ with $x-y=nz$ for some $z\in A$ one has
$f(x)=f(y)$.  
For $B\subset A$ and $f\in {\cF}(A,S)$
we denote by $f_B$ the restriction to $B$ 
(or simply $f$ if the domain 
is clear from the context). 

\begin{assu}\label{assu:a}
Throughout, all abelian groups are 
either in ${\mathcal A}_0$ or in ${\mathcal A}_q$. 
\end{assu}

\subsection{Definitions}
\label{sect:defn}

We will work in the following setup:
$$
A\subset V\subset K
$$
where $A$ is a  
$\Z$-(resp. ${\bf F}_q$-) sublattice of a  
$k$-vector space $V$ which is embedded into $K$.

\begin{defn}\label{defn:inv}
Let $A\in {\mathcal A}_0$ 
(resp. $A\in {\mathcal A}_q$) 
and  $f\in \cF(A,S)$. We say that $f$ is invariant if
$$
f(na)=f(a)
$$ 
for all $a\in A$ and all 
$n\in \Z_{(0)}$ (resp. all $n\in \Z_{(q)}$)

Let $V$ be a vector space over $k$ and $f\in \cF(V,S)$. 
We say that $f$ is invariant if
$$
f(\kappa v)=f(v)
$$
for all $\kappa\in k^*$ and $v\in V$. 
\end{defn}

An invariant function on $A={\bf F}_q^n$ (minus $0_A$)
can be considered 
as a function on $\P^{n-1}({\bf F}_q)= (A\setminus 0)/{\bf F}_q^*$. 
An invariant function on $V\setminus 0$ can be considered as
a function on the projective space $\P(V)$ (over $k$) and we will 
denote by $\cF(\P(V),S)\subset \cF(V,S)$ the space of such functions.

\begin{defn}[Filtration]
\label{defn:filtered} 
Let $A$ be finitely generated and $\cI$ a totally ordered set. 
A (strict) {\em filtration} on $A$ with respect to $\cI$   
is a set of subgroups 
$A_{\iota}\subset A$ (with ${\iota\in \cI}$) 
such that 

$\bullet$ $A=\cup_{\iota\in \cI}A_{\iota}$;

$\bullet$ if $\iota <\iota'$ then
$A_{\iota'}$ is a proper subgroup of $A_{\iota}$.
\end{defn}

\begin{nota}
\label{nota:ol}
Denote by 
$$
\overline{A}_\iota:= 
A_\iota\setminus \cup_{\iota' >\iota} A_{\iota'}.
$$
\end{nota}

Notice that for all $\iota\in \cI$ we have 
$ \overline{A}_\iota\neq \emptyset.$

\begin{defn}[AF-functions]
\label{defn:af} 
A function $f$ on a finitely generated group 
$A$ (as in \ref{assu:a}) is 
called an {\em abelian flag} function if 

$\bullet$\,\, $f$ is invariant;

$\bullet$\,\, $A$ has a (strict) filtration by groups
(with respect to an ordered set $\cI$) such that 
$f$ is constant on $\overline{A}_\iota$
for all $\iota \in \cI$.

If  $A$ is not finitely generated then     
$f$ is an abelian flag function if $f_{A'}$ is an abelian 
flag function for every finitely generated subgroup $A'\subset A$.  
\end{defn}

We denote the set of abelian flag functions by $\AF(A,S)$. 
This property does not depend on the value of $f$ on 
the neutral element $0_A$. We will identify functions
which differ only on $0_A$.

\begin{defn}\label{defn:gf}
Let $V$ be a vector space over $k$. A function $f\in \cF(V,S)$ is 
called an abelian flag function if 

$\bullet $ $f$ is invariant;

$\bullet$ for all (additive) sublattices $A\subset V$
the restriction  $f_A\in \AF(A,S)$.
\end{defn}

\begin{defn}[c-pairs]
\label{defn:c-pair}
Let $A$ be an abelian group as above and $S$ a ring.
We will say that $f_1,f_2\in \cF(A,S)$ form a c-pair
if for every subgroup $C\subset A$ of rank 2 
one has
$$
\rk(\langle f_1,f_2,1\rangle_S)\le 2.
$$ 
\end{defn}

\begin{defn}[LF-functions]
Let $V$ be a unital algebra over $k$ and $S$ an abelian group. 
A function $f\in \cF(V,S)$ is called a {\em logarithmic} function
if 

$\bullet$ $f$ is invariant;

$\bullet$  $f(v\cdot v')=f(v)+f(v')$ for all $v,v'\in V\setminus 0$.   

\end{defn}

The set of logarithmic functions will be denoted by $\LF(V,S)$. 
We shall refer to abelian flag  (resp. logarithmic) functions
as AF-functions  (resp. LF-functions).

\subsection{First properties}
\label{sect:1prop}

\begin{rem}\label{rem:finitely-generated}
Assume that $A$ is {\em finitely generated}. 
Then for every  $f\in \AF(A,S)$
and every subgroup $B\subset A$ 
there exists a proper 
subgroup $B_f^1\subset B_f^0=B$
such that $f$ is constant on 
the complement $B\setminus B_f^1$. 
In particular, 
if $b_0\in B\setminus B_f^1$ and $b_1\in B_f^1$ then 
$f(b_0+b_1)=f(b_0)$. 
Thus we can speak about generic elements of $B$ and  
the generic value of $f$ on $B$. 
We obtain a decreasing (possibly finite) 
$\N$-filtration $(A_f^n)$ on $A$: 
$A_f^n$ is the {\em subgroup} of nongeneric
elements in $A_f^{n-1}$. 
Notice that an analogous statement
for infinitely generated groups 
is not true, in general (for example, 
$\Q$ and valuation subgroup for a nonarchimedean valuation).   
\end{rem}

\begin{lemm}\label{lemm:easy-p}
Assume that $f\in \AF(A,S)$. Then 

$\bullet$ for all subgroups  $B\subset A$ 
one has $f_B\in \AF(B,S)$; 

$\bullet$ if $S$ is a ring then 
$sf+s'\in \AF(A,S)$ for all $s,s'\in S$;  

$\bullet$ for every map $h\,:\,S\to S'$ 
one has $h\circ f\in \AF(A,S')$. 
\end{lemm}

\begin{proof}
Evident from the definition. 
\end{proof}

\begin{lemm}
\label{lemm:index}
Let $A\in {\mathcal A}_0$ and 
$f\in \cF(A,S)$. Let $B\subset A$
be a subgroup of finite index. 
If $f_B\in \AF(B,S)$ then $f\in \AF(A,S)$. 
\end{lemm}

\begin{proof}
Observe that $nA\subset B$ for some $n\in \N$. 
By Lemma~\ref{lemm:easy-p},
$f_{nA}\in \AF(nA,S)$. Now use the invariance of $f$.
\end{proof}

\subsection{Orders}
\label{sect:orders}

\begin{rem}\label{rem:partial}
Let $f\in \cF(A,S)$ be such that
for every subgroup $B\subset A$ with $\rk(B)\le 2$ 
the restriction $f_B\in \AF(B,S)$.
Then $f$ defines a partial relation $\tilde{>}_f$ 
on $A$ as follows: let $b,b'\in A$ with $f(b)\neq f(b')$ 
and consider the subgroup $B=\langle b,b'\rangle$. 
One of these elements, say $b$, 
is generic in $B$ and the other is not. Then we  
define $b \tilde{>}_f b'$.
\end{rem}

\begin{lemm}[Definition]\label{lemm:part-ord}
A function  $f\in \AF(A,S)$ defines an {\em order} 
$>_{f,A}$ on $A$ as follows:

$\bullet$\,\, 
if $f(a)\neq f(a')$  then $a>_{f,A} a'$ iff 
$a\tilde{>}_f a'$ (that is, $f(a+a')=f(a)$);

$\bullet$\,\,
if $f(a)= f(a')$ then $a>_{f,A} a'$ iff there exists a $b\in A$ 
(a separator) such that $f(a)\neq f(b)$ and 
$a>_{f,A} b >_{f,A} a'$; 

$\bullet$\,\, 
finally, $a =_f a'$ if for all $a'' $ we have
$a''>_{f,A} a$ iff $a''>_{f,A} a'$. 
\end{lemm}

\begin{nota}
\label{nota:rel}
We will write $>_f$ or $>$ and $\tilde{>}$ 
whenever $A$ and $f$ are clear
from the context. We will also use the symbols $\ge$ and $\ge_f$.  
\end{nota}

\begin{proof}
If $\rk(A)<\infty$, then for every element $a\neq 0$
there exists an $n(a)$ 
such that $a\in A_f^{n(a)}\setminus A_f^{n(a)+1}$. 
Then $a> b$ iff $n(a)<n(b)$.  
For general $A$, the 
correctness of the definition and the transitivity of
$>$ are checked on finitely generated subgroups of $A$.
For correctness, 
we assume that $f(a)=f(b)$ and consider the possibility
that 
$$
a\, \tilde{>}\, c\, \tilde{>}\, b\, \tilde{>}\, c'\,\tilde{>}\, a
$$ 
for two separators $c,c'$
(leading to a contradiction). For transitivity, we
may need to consider the possibility 
$$
a > b> c\ge_f a
$$
(with separators, if necessary). 
\end{proof}

Let $A$ be an abelian group and $f\in \AF(A,S)$. 
Define the subsets:
$$
A_f^{\alpha}= \{ a\in A\,|\, a\ge_f \alpha\}
$$
where $\alpha$ is (a representative of) 
the equivalence class in 
$A$ with respect to the equivalence relation $=_f$.

\begin{lemm}\label{lemm:ge-filtration}
Assume that $\rk(A)<\infty$ and $f\in \AF(A,S)$. 
Then each $A_f^{\alpha}$ is a subgroup of $A$ and
$(A_f^{\alpha})$ 
is a filtration on $A$ in the sense of 
Definition~\ref{defn:filtered}.
Moreover, $f$ is constant
on $\ovl{A}_f^{\alpha}$ for all $\alpha$.
\end{lemm}

\begin{proof} 
Evident. 
\end{proof}

\begin{rem} 
If $\rk(A)<\infty$ then $(A_f^\al)$ coincides
with the filtration $(A_f^n)$
introduced in Remark~\ref{rem:finitely-generated}.
Notice that the filtration $(A_f^\al)$ is not functorial under
restrictions.  
In general, if $f\in \AF(A,S)$ and $B\subset A$ 
is a proper subgroup then 
$A_f^\al\cap B\neq B_f^\al$. 
In general, for maps $h\,:\, S\ra S'$ the filtration 
$A_f^\al$ does not coincide with 
$A_{h\circ f}^\al$. However, $A_{h\circ f}^\al$ can
be reconstructed starting from $A_f^\al$.
We will be interested in the case when
$S'=\Z/2$. 
\end{rem}

\begin{lemm}\label{lemm:part}
Let $f\in \cF(A,S)$ be such that:

$\bullet$ for every $B$ with $\rk(B)\le 2$ the restriction
$f_B\in \AF(B,S)$ (by Remark~\ref{rem:partial},
this defines a partial relation $\tilde{>}$ on $A$);

$\bullet$ $\tilde{>}$ extends to an order 
$>$ on $A\setminus 0$ (transitivity);

$\bullet$ for all $a,a',a''\in A$ such that
$a> a'$ and $a> a''$ one has 
$$
a> a'+a''.
$$

Then $f\in \AF(A,S)$.  
\end{lemm}

\begin{proof}
Evident. As in Lemma~\ref{lemm:ge-filtration}, 
we obtain a filtration by groups. 
\end{proof}

\subsection{Rank 2 case}
\label{sect:first}

\begin{exam}\label{exam:typical}
A typical AF-function is given as follows.
Let $p$ be a prime number and $A=\langle a,a'\rangle$ with 
$$
\begin{array}{ccl}
A_f^{2n}    & = &   \Z p^na \oplus    \Z p^{n}a' \\ 
A_f^{2n+1}  & = &  \Z p^n a \oplus \Z  p^{n+1}a'
\end{array}
$$
with $f$ taking two values on $A$: one value on 
$\overline{A}_f^{2n}$ and a different value on  
$\overline{A}_f^{2n+1}$.
\end{exam}

\begin{lemm}\label{lemm:z2}
Let $A=\Z\oplus \Z$ and $f\in \AF(A,S)$. Then 
$f$ is one of the following

$\bullet$ $f$ is constant on $A\setminus 0$;

$\bullet$ $f$ is constant on $A\setminus \Z a$, for some $a\in A$;

$\bullet$ there exists a prime number $p$ and a 
subgroup $C$ of index $p^k$
(for some $k\ge 0$) such that $f$ is constant on $A\setminus C$
and $f_C$ is as in the Example~\ref{exam:typical}.  
\end{lemm}

\begin{nota}
\label{nota:p}
In the second  (resp. the third) case we put  
$$
p(A)=p(A,f):=0,
$$ 
$$
p(A)=p(A,f):=p.
$$ 
\end{nota}

\begin{proof}
Assume that $f$ is nonconstant on $A\setminus 0$. 
Then there exist two primitive elements
$a,a' \in A$ such that $f(a)\neq f(a')$ 
(and $B:=\langle a,a'\rangle $ is a subgroup of finite index in $A$). 
Then one of these generators, say $a'$, lies in 
the subgroup $B_f^1$. 
This means that $B/B^1_f$ is a cyclic group. 
If it is a free cyclic group, 
then the function $f$ is of the second type.

If it is a finite group then  
there is a proper subgroup $C$ of finite index in $B$
such that $C$ contains $B^1_f$ and $[C:B^1_f]=p$ (for some prime $p$). 
The function $f$ is constant on $C\setminus B^1_f$. 
Hence $C^1_f=B^1_f$. 
We have the diagram
$$
\begin{array}{cccccccclc}
A &\supset & B &\supset & C &\supset & B^1_f  & \supset  & B^2_f    & \supset ... \\
 &        & &        & &        &  ||    &          & ||       &   \\
 &        & &        & C &\supset &  C^1_f & \supset  & C^2_f=pC &     
\end{array}
$$
Indeed, since the generic value of $f$ on 
$pC$ is equal to the generic value of $f$ on $C$
it is not equal to the generic value of $f$ on $C^1_f$. 
It follows that
$pC\subset C^2_f$ and since $[C^1_f : pC]=p$ 
we have $pC=C^2_f$.
By invariance, the function 
$f_C$ is as in the Example~\ref{exam:typical}. 
Again by invariance, 
the index $[A:C]$ is a $p$-power and $f$ (on $A$) is of the
third type.

\end{proof}

\begin{coro}
\label{coro:rk2-easy}
Let $A$ and $f$ be as in Lemma~\ref{lemm:z2}. Then 

$\bullet$ $f$ takes at most two values on $A\setminus 0$;

$\bullet$ if $na\in A^1_f$ for some $n$ with 
$\gcd(n,p(A))=1$ then $a\in A^1_f$; 

$\bullet$ if $A=\langle a_1,a_2\rangle$ and $a_1-a_2=p(A)a_3$
($a_3\in A$) then $f(a_1)=f(a_2)$;

$\bullet$ if $|A/A^1_f|=2$ and 
$a\in A$ is primitive then $f(a+2a')=f(a)$ for all $a'\in A$;

$\bullet$ if $|A/A^1_f|>2$ then
$A$ has a basis $\{a_1,a_2\}$ such that
all three elements $a_1,a_2,a_1+a_2$ are generic. 
\end{coro}

\begin{lemm}\label{lemm:z2-p}
Let $A=\Z/q\oplus \Z/q$ (with $q$ prime) 
and $f\in \AF(A,S)$. 
Then $f$ (considered as a function on $\P^1({\bf F}_q)$) 
is constant on the complement to 
some point $P\in \P^1({\bf F}_q)$. 
\end{lemm}

\begin{proof}
See the proof of Lemma~\ref{lemm:z2}. 
\end{proof}

\begin{lemm}
\label{lemm:agf}
Let $A$ and $f\in \cF(A,\Z/2)$ be such that:

$\bullet$ $A$ has a basis $(a,b)$ with 
$f(a)=f(a+b)\neq f(b)$;

$\bullet$ $f$ is invariant (cf. Definition~\ref{defn:af});

$\bullet$ $f$ satisfies a 
functional equation: 
for all $a',b'\in (A\setminus 0)$ 
with 
$$
f(a')=f(a), \,\,f(b')=f(b)\,\,\, {\rm  and}\,\,\, f(a')=f(a'+b')
$$ 
one has
\begin{equation}
\label{eqn:fe}
f(ma'+nb')=f(ma'+(n+km)b')
\end{equation}
for all $k,m,n\in \Z$.

Then $f\in \AF(A,\Z/2)$. 
\end{lemm}

\begin{proof} 
We have a decomposition $A=A_a\cup A_b$ 
into two subsets (preimages of $0,1$).  
We will generally use the letter $a$ for elements in $A_a$
and $b$ for elements in $A_b$.
Thus $f(a')=f(a)\neq f(b)=f(b')$ for all $a,a',b,b'\in A$.

First we consider the case $A\in \cA_q$.
By the functional equation,
$$
f(a+nb)=f(a)
$$ 
for all $n\in \Z$. 
By invariance, 
$$
f(ma+nb)=f(a)
$$ for all $m,n\in \Z_{(q)}$. 
Thus $f$ is constant on $A\setminus A^1_f$
(where $A^1_f=\Z b$). 

Now we turn to the case $A\in \cA_{0}$. 
Denote by $A_f^1=\langle A_b\rangle \subset A$ 
the subgroup in $A$ generated by elements 
$b'\in A_b$.
We claim that $A_f^1$ is 
a {\em proper} subgroup of $A$
(and clearly, $f(a_1)=f(a)$ 
for all $a_1\in (A\setminus A^1_f)$).

Consider a pair of elements 
$$
\begin{array}{ccc}
b_1& = & m_1a+n_1 b, \\ 
b_2& = & m_2a+n_2b.
\end{array}
$$ 
We can assume that $m_1,m_2>0$. 
Let 
$$
d_A:=\min (\gcd(m_1,m_2))
$$ 
be the minimum over all pairs 
$(b_1,b_2)\in A_b\times A_b$ (with positive $m_1,m_2$).

Assume first that
there exists a $b_1$ such that 
$b_1=d_A a+n_1b$ (for some 
$n_1\in \Z$). 
This is impossible for $d_A=1$, by the functional equation. 
Now consider the case $d_A>1$. 
In this case $A^1_f$ is a proper subgroup 
of $A$ (since for all $b_2=m_2a+n_2b$ 
the coefficient $m_2$ is divisible by $d_A$ and 
consequently for all $ma+nb \in A^1_f$
the coefficient $m$ is divisible by $d_A$).

Now assume that there are no such $b_1$. 
Choose a pair 
$$
\begin{array}{ccc}
b_1 & = & d_Am_1a+n_1b,\\ 
b_2 & = & d_Am_2a+n_2b
\end{array}
$$ 
as above (such that  
$(m_1,m_2)=1$) 
and integers $l_1,l_2\in \Z$ 
such that $m_1l_1+1=m_2l_2$.
Then, (using invariance), 
\begin{align*}
f(d_A r a + e_1b) & = f(d_A(r+1)a+e_2b)\\
                  & = f(b)
\end{align*}
for some $e_1,e_2\in \Z$. 
Pick the smallest positive $r_0$ with this property.
Then $r_0>1$ (since 
$f(d_Aa+nb)=f(a)$ for all $n$, by assumption)
and 
$$
f(d_A(r_0-1)a + nb)=f(a)
$$ 
for all $n\in \Z$. 
Therefore, (using functional equation and invariance),  
\begin{align*}
f(a) & = f(d_A(r_0-1)a + (2e_1-e_2)b) \\
     & = f((-d_A a + (e_1-e_2)b) + (d_Ar_0a +  e_1b))\\
     & = f(d_A(r_0+1)a + e_2b)\\
     & = f(b),
\end{align*}
contradiction.

Thus $A^1_f$ is a proper subgroup 
of $A$ and $a_1\in A_a$ for all $a_1\notin A^1_f$. 
Now consider the subgroup 
$$
A^2_f:=\langle A^1_f\cap A_a\rangle 
$$ 
(generated by $a_2\in A^1_f$ with $a_2\in A_a$). 
We claim that $A^2_f$ is a proper subgroup of $A^1_f$. 
(Warning: the conditions of the 
Lemma are {\em not} symmetrical
with respect to $a$ and $b$.) 

First observe that there exists 
a basis $a',b'$ of $A^1_f$ such that 
$$
f(a'+ b')=f(b)
$$
(otherwise, the subgroup $\langle A_b\rangle$ would 
be a proper subgroup in $A^1_f$, contradiction). 
We fix such a basis and claim that 
$$
f(ra' + b')=f(b)
$$ 
for all $r\in \Z$. 
Indeed, if this is not the case then
there exists a positive integer $r_0>1$ such that
$$
\begin{array}{ccc}
f(a)           & = & f(r_0a'+b'), \\
f((r_0-1)a'+b')& = & f(b).
\end{array}
$$  
Then 
\begin{align*}
f(a) & = f(a' + ((r_0-1)a' + b'))\\
     & = f(a'-((r_0-1)a' + b')) \\
     & = f((r_0-2)a'+b')
\end{align*}  
(here we used the functional equation, 
and the invariance). This contradicts the 
minimality of $r_0$. Similar argument works for $r<0$.

Consider the set of pairs
$$
\begin{array}{ccc}
a_1 & = & m_1a'+n_1b' \\
a_2 & = & m_2a'+n_2b'
\end{array}
$$ 
with positive $n_1,n_2$ and denote by 
$d_B$ the smallest $\gcd(n_1,n_2)$ on this set.

Assume that $d_B>1$ and that there exists an 
$a_1=m_1a'+d_Bb'$.
Then for every $a_2=m_2a'+n_2b'\in A^1_f$
the coefficient $n_2$ is divisible by $d_B$ and  
$A^2_f$ is a proper subgroup.

Now we can assume that 
$$
f(ma'+d_Bb')=f(b)
$$ 
for all $m\in \Z$.
Choose a pair $a_1, a_2$ such that 
$\gcd(n_1,n_2)=d_B$. 
Then, (by invariance), 
\begin{align*}
f(a) & = f(r_1 a' + e d_B b') \\
     & = f(r_2a'+(e+1)d_Bb')
\end{align*}
for some $r_1,r_2\in \Z$ and $e>0$. 
Pick the smallest $e_0>1$ with this property. 
Then, (using 
the fact that 
$$
f((r_2-r_1)a' +d_B b')=f(b)
$$ 
and the functional equation), we get
\begin{align*}
f(a) & = f((r_1 a' + e_0 d_B b') +(r_2 -r_1)a' + d_Bb') \\
     & = f((2r_1-r_2)a'+(e_0-1)d_Bb').
\end{align*}
This contradicts the minimality of $e_0$. 

It follows that the subgroup $A^2_f\subset A^1_f$ is
a proper subgroup and $f$ takes the value $f(b)$ on 
the complement $A^1_f\setminus A^2_f$. 

Since $A^2_f$ has a basis $(a'',b'')$
with 
$$
f(a'')=f(a''+b'')\neq f(b''),
$$ 
we can apply the inductive step to $A^2_f$.
\end{proof}

\subsection{Rank 3 case: $\cA_q$}
\label{sect:rk3}

\begin{prop}
\label{lemm:red2-p}
Let $q>2$ be a prime number, 
$A=\Z/q\oplus \Z/q\oplus \Z/q$ and 
$f\in {\mathcal F}(A,\Z/2)$. Assume that
for all subgroups $C\subset A$ 
with $\rk(C)\le 2$ we have $f_C\in\AF(C,\Z/2)$. 
Then $f\in\AF(A,\Z/2)$.
\end{prop}

\begin{proof}
We can consider $f$ as a function on $\P(A)$.  
By Lemma~\ref{lemm:z2-p}, for every $C\subset A$ with
$\rk(C)=2$ the restriction $f_C$ is either constant on 
$\P(C)\subset \P(A)$ or constant everywhere except one point.
Let $L_i$ be the set of lines $\ell\subset \P(A)$ such that
the generic value of $f_{\ell}$ is $i$ (for $i=0,1$).  
Assume that $f$ is nonconstant on $\P(A)$. 
If $L_0$ is empty then there exists only one point $P\in \P(A)$ 
with $f(P)=0$ (otherwise we can draw a line of type $L_0$ through
two such points, and $0$ must be the generic
value on this line, contradiction). 
In this case $f\in \AF(A,\Z/2)$.   
Thus we can assume that both $L_0$ and $L_1$ are nonempty and 
that, for example, $|L_0|\ge |L_1|$. Then  $|L_0|\ge (q^2+q+1)/2$. 
Choose an $\ell_1\in L_1$. 

There are two cases: $f$ is constant on $\ell_1\in L_1$
or $f$ is nonconstant on some line $\ell_1$.
In both cases there exists at least one point 
$P\in \ell_1$ such that $f(P)=1$ and such that there are two 
distinct lines 
$\ell_0,\ell_0'$ passing through $P\subset \ell_1$. 

Indeed, assume in both cases
that through every generic point of $\ell_1$ there
passes only one line of type $L_0$. In the first case  
the total number of lines of type $L_0$ is bounded by
$q+1$, contradiction to the assumption that $|L_0|\ge (q^2+q+1)/2$.
In the second case, there are at most 
$q$ lines of type $L_0$ passing
through the nongeneric point and, by assumption,
at most $1$ line of type $L_0$ 
passing through each of the remaining $q$ 
generic points of $\ell_1$ 
(every line in $L_0$ intersects $\ell_1$ in one point).
Thus their number is bounded by $2q<(q^2+q+1)/2$, contradiction.

For any pair of points $Q\in\ell_0 \setminus P$, $Q'\in
\ell_0'\setminus P$ we have $f(Q)=f(Q')=0$. 
The lines through $Q,Q'$ are all of type $L_0$. 
Pick a point $P'\in \ell_1$ such that $P\neq P'$ and 
$f(P')=1$ (such a point exists since $\ell_1$ has 
at least 3 points and the generic value of $f$ on $\ell_1$ is 1). 
Every line through $P'$ which does not pass through $P$ is 
of type $L_0$ (since it intersects $\ell_0,\ell_0'$ 
in distinct points). 
The family of such lines covers $\P(A)\setminus \ell_1$. 
It follows that the value of $f$ on $\P(A)\setminus \ell_1$ is 0
and that $f\in \AF(A,\Z/2)$.
\end{proof}

\subsection{Exceptional lattices}
\label{sect:exc}

\begin{exam}
\label{exam:exam}
Let $\bar{A}=\Z/2\oplus\Z/2\oplus \Z/2$ and $\bar{f}\in \cF(\bar{A},\Z/2)$.
Assume that for all subgroups $\bar{C}\subset \bar{A}$ of rank 2 one
has $f_{\bar{C}}\in \AF(\bar{C},\Z/2)$ but 
$\bar{f}\notin \AF(A,\Z/2)$. Then $A$ has a basis 
$\bar{e}_1,\bar{e}_2,\bar{e}_3$ such that 
$$
f(\bar{e}_1+ \bar{e}_3)=f(\bar{e}_2+\bar{e}_3)=
f(\bar{e}_1+\bar{e}_2+\bar{e}_3)=0\neq f(\bar{x})
$$
for all other $\bar{x}$ (up to addition of 1 modulo 2). 

Indeed, 
since $\P^2(\Z/2)$ has seven points it suffices to assume that $f$ takes
the same value on three of them and a different value 
on the remaining four. 
If the three points are on a line we have an AF-function. 
If not we get the claim. (In particular, such an $f$ contradicts
the conclusion of Lemma~\ref{lemm:red2-p}.)

\vskip 0,5cm
\centerline{\includegraphics{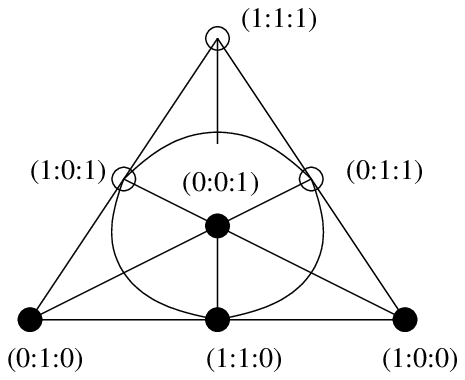}}
\vskip 0,5cm

Any function on $A=\Z\oplus \Z\oplus \Z$ induced 
from the function on $\bar{A}$ considered above
has the property that for every $C\subset A$ of rank $\le 2$
$f_C\in \AF(C,\Z/2)$.
\end{exam}

We give another example of a function on $\Z^3$ 
with the same property.

\begin{exam}\label{exam:exam2}
We keep the notations of Example~\ref{exam:exam}.
Choose a basis ${e}_1,{e}_2,{e}_3$ of $A=\Z^3$. 
Consider the projection $A\ra \bar{A}= A/2A$, 
taking ${e}_j$ to $\bar{e}_j$.
The function $f$ is defined by its values on primitive elements
$$
a=n_1e_1+n_2e_2+e_3n_3.
$$ 
If 
$(n_1,n_2,n_3)\neq (0,1,0)$ modulo 2
then 
$$
f(a)=\bar{f}(\bar{a})
$$
(where $\bar{f}$ was defined in Example~\ref{exam:exam}).
Otherwise, 
$$
f(n_1e_1+n_2e_2+n_3e_3)=\left\{ 
\begin{array}{ccc} 0 & {\rm if} & n_1=0 \mod 4\\
                   1 & {\rm if} & n_1=2 \mod 4
                                \end{array} \right.
$$

\vskip 0,5cm
\centerline{\includegraphics{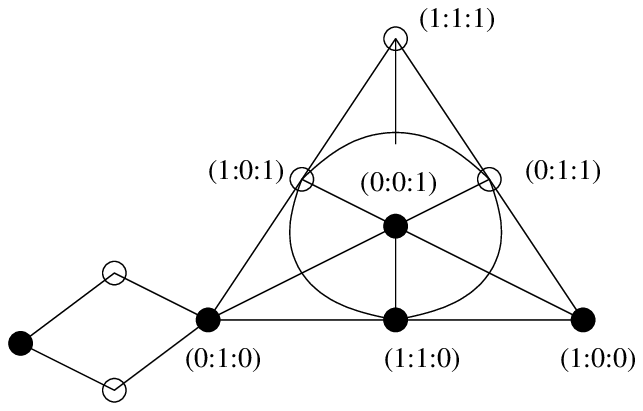}}
\vskip 0,5cm

\end{exam}

\subsection{Rank 3 case: $\cA_0$}

\

Let $A$ be an abelian group and $f\in \cF(A,\Z/2)$. 
We have a decomposition of the set 
$A=A_a\cup A_b$ (preimages of 0 or 1, respectively).
We will say that $A$ has a {\em special basis}
(with respect to $f$) if 
$A=\langle a_1,a_2,b_1\rangle$ with 
$$
a_1,a_2,a_1+b_1, a_2+b_1\in A_a, \,\,\, b_1\in A_b.
$$

\begin{prop}\label{prop:red2-0}
Let $A=\Z\oplus \Z\oplus \Z$ and 
$f\in {\mathcal F}(A,\Z/2)$. 
Assume that $A$ has a special basis 
(with respect to $f$) and that 
for every  subgroup $C$ with $\rk(C)\le 2$ one has
$f_C\in \AF(C,\Z/2)$.
Then there is a proper subgroup $A^1\subset A$ such 
that $f$ is constant on $A\setminus A^1$. 
\end{prop}

\begin{proof}
The proof is subdivided into a sequence of Lemmas.

\begin{lemm}
\label{sublemm:basis}
Assume $A$ has a special basis $\{ a_1,a_2,b_1\}$. 
Then $A$ does not have a basis $\{ b_1, b_2, b_3 \}$ with
$b_1,b_2,b_3\in A_b$.  
\end{lemm}

\begin{proof}
Assume the contrary. Then 
$$
C:=\langle a_1+e_1b_1,a_2+e_1'b_1\rangle 
= \langle b_2,b_3 \rangle
$$
for some $e_1,e_2'\in \Z$. 
We know that 
$$
f(a_1+e_1b_1)=f(a_2+ e_1'b_1)=f(a)
$$ 
for all $e_1,e_1'$.
Contradiction to, for example, \ref{lemm:z2}:
$C$ cannot have such pairs of generators.   
\end{proof}

\begin{lemm}
Assume that $A$ has a special basis $\{ a_1,a_2,b_1\}$. 
Then  
$$
\langle A_b\rangle \subset A
$$  
is a proper subgroup. 
\end{lemm}

\begin{proof}
Consider the projection
$$
{\rm pr}\,:\, A\ra \hat{A}:=A/\langle b_1\rangle.
$$
Assume that there exists an element 
$$
b=n_1a_1+n_2a_2 + m_1b_1\in A_b
$$
such that $\hat{b}=n_1a_1+n_2a_2$ is primitive in 
$\hat{A}=\langle a_1,a_2\rangle$. 
Then it is part of a basis 
$\{ \hat{x},\hat{b}\}$ of $\hat{A}$. 
Take any $x$ in the preimage ${\rm pr}^{-1}(\hat{x})$. 
Then 
$$
A=\langle x,b,b_1\rangle.
$$ 
By Lemma~\ref{sublemm:basis}, 
$x\notin A_b$, so we will denote it by $a$. 
Assume that 
$$
a+mb+m_1b_1\in A_b
$$ 
for some $m,m_1\in \Z$. 
This contradicts \ref{sublemm:basis},
since 
$$
A=\langle a+mb+m_1b_1,b,b_1\rangle.
$$.

Consider the set $R$ of all $r\in \N$ such that 
$$
b_r:=ra+ mb+m_1b_1\in A_b
$$
(for some $m,m_1\in \Z$). We have seen that $r>1$. 
We claim that $r,r'\in R$ implies that $g:=\gcd(r,r')\in R$.
Indeed, assume the contrary and 
choose $l,l'\in \Z$ so that $g=lr-l'r'$. 
By invariance,
$lb_r, l'b_{r'}\in A_b$. 
In the subgroup
$B:=\langle b_r,b_{r'}\rangle$ the element 
$a_g:=lb_r-l'b_{r'}\in A_a\cap B$ is nongeneric.  
This implies that 
$$
b_r-na_g \in A_b\cap B
$$
for all $n$. This leads to a contradiction 
and the claim follows. 
Thus we have proved that for all 
$b'\in A_b$ the corresponding 
coefficients $n_1'$ are either zero 
or have a common divisor $>1$. 
Consequently, $\langle A_b\rangle$ is a proper subgroup.

Now we assume that ${\rm pr}(A_b)$ does not
contain primitive elements of $\langle a_1,a_2\rangle$, 
in other words: for all 
primitive $a\in \langle a_1,a_2\rangle$ and all $m_1\in \Z$
one has 
$$
a+m_1b_1\in A_a.
$$ 
For two pairs of $(a,m_1)$ and $(a',m'_1)$ 
with primitive $a,a'$ such that 
$$
\langle a_1,a_2\rangle =\langle a,a'\rangle
$$
consider the subgroups
$$
\begin{array}{ccc}
D        & := & \langle a,m_1b_1\rangle\\
D'       & := & \langle a',m'_1b_1\rangle
\end{array}
$$
and assume that both $p=p(D),p'=p(D')\neq 0$.

We claim that $p=p'$. Indeed, assume the contrary. 
By Lemma~\ref{lemm:z2}, there exist integers $k,k'$ such that
$$
f(q a+ m_1b_1)=f(q' a'+ m'_1b_1)=f(b_1),
$$ 
where $q=p^k,q'={p'}^{k'}$ for some $k,k'\in \N$. 
Now consider the group 
$$
E:=\langle qa+m_1b_1, q'a'+m'_1b_1\rangle.
$$ 
For all $n_1\in \Z$ coprime to $p'$ and all $n'_1\in \Z$ 
coprime to $p$ the element
$$
n_1(qa+m_1b_1)+n'_1(q'a'+m'_1b_1)=
n_1qa+ n'_1q'a' + 
(n_1m_1+n'_1m'_1)b_1 \in A_a,
$$
(since ${\rm pr}(A_b)$ does not contain
primitive elements).  
The subset of such elements cannot be contained in 
a proper subgroup of $E$. On the other hand, 
it has to be: both generators of $E$ are in $A_b$
and $f\in \AF(E,\Z/2)$. Contradiction to the assumption that
$p\neq p'$. 
Since for any pair of primitive $a,a'$ generating
a sublattice of finite index in $\langle a_1,a_2\rangle$ there
exists a primitive element $a_0$ such that 
$$
\langle a,a_0\rangle = \langle a',a_0\rangle=\langle
a_1,a_2\rangle
$$
we conclude that for the corresponding $D$ 
as above either $p(D)=0$ or
$p(D)=p$ for some fixed prime $p$.

To finish the proof of the lemma, 
consider an element
$$
na +m_1b_1\in A_b
$$ 
for some $n>1$, some
primitive $a\in \langle a_1,a_2\rangle$ and some $m_1$
(coprime to $n>1$). 
There are two possibilities: 
either $n$ is zero or $n$ is divisible by
a fixed prime $p$ (which is independent of the coefficients). 
It follows that 
$\langle A_b\rangle\subset A$ is a proper subgroup.
\end{proof}

This concludes the proof of Proposition~\ref{prop:red2-0}
\end{proof}

\begin{prop}\label{prop:special}
Let  $A=\Z\oplus \Z\oplus \Z$ and $f\in \cF(A,\Z/2)$.
Assume that $A$ does not have a special basis 
(with respect to $f$)
and that for all subgroups $C\subset A$ of rank 2
one has $f_C\in \AF(C,\Z/2)$. 
Then $f$ (up to addition of 1 modulo 2) is 

Case 1:  as in Example~\ref{exam:exam} or 

Case 2:  as in Example~\ref{exam:exam2}.
\end{prop}

\begin{proof}
Assume first that for all subgroups $C\subset A$ of rank 2
the function $f$ is either constant on $C\setminus 0$ or
$|C/C^1_f|=2$. Then $f$ is induced from $A/2A$. Indeed, 
consider the subgroup 
$C:=\langle x,y\rangle$ (with $x\in A$ primitive and $y$
nonproportional to $x$). It suffices to consider
the case when $f$ is nonconstant and thus 
$|C/C^1_f|=2$. By  Corollary~\ref{coro:rk2-easy},
if $x-y=2z$ we have
$$
f(x)=f(y),
$$ 
so that $f$ is induced from $A/2A$. 

Now we assume that there exists a subgroup 
$C$ of rank 2 such that $|C/C^1_f|>2$ and $f$ has generic value, 
say $f(a)$ on $C$. 
By Corollary~\ref{coro:rk2-easy}, we can choose a basis 
$C=\langle a_1,a_2\rangle$ such that all three 
$a_1,a_2,a_1+a_2\in A_a$.
We will fix such a basis. 
For any $d\in A$ such that $\langle d, C \rangle = A$ 
consider the shift $d+C\subset A$. 

\

1. We claim that both
$$
(d+C)\cap A_a\neq \emptyset \,\,\, {\rm and}\,\,\,
(d+C)\cap A_b\neq \emptyset.
$$
Indeed, assume this is not so and 
choose, in the first case, 
some element $b\in d+C$. Then $\{ b +a_1, b+a_2, a_1\}$ 
is a special basis of $A$. In the second case, for 
any $a\in d+C$ we get a special basis $\{ a, a_1, b_1\}$, 
(where $b_1\in A_b$ is some nongeneric element in $C$).
Contradiction. 

\

For any pair of generators $\{ a_1',a_2'\} $ of $C$
(without the assumption that $a_1'+a_2'\in A_a$)
we have:

\

2. (Forbidden triangle.) There are no $b\in (d+C)\cap A_b$ 
such that both 
$$
b+a_1',b+a_2'\in A_a.
$$   
Indeed, $\{ a_1',a_2',b\} $ would be a special basis for $A$. 

\

3. (Forbidden square.) There are no $b\in (d+C)\cap A_b$
such that all three 
$$
b+a_1',b+a_2',b+a_1'+a_2'\in A_b.
$$
Indeed, $\{ b, b+a_1', a_2'\}$ would be a special basis for $A$. 

\

4. Choose any element $b\in d+C$ and consider the subset
$\{b+na_1\}$ (with $n\in \N$). Than there are two possibilities:
either $b+na_1\in A_b$ for all $n\in \N$ 
or there exists an $n_0>0$ such that 
$$
b+n_0a_1\in A_b\,\,\, {\rm  and}\,\,\, b+(n_0+1)a_1\in A_a.
$$  
Let us consider the second case:
rename $b+n_0a_1$ to $b$. 

Then $b+a_1\in A_a$. 
By 2, $b+a_2,b-a_2\in A_b$. By 2 and by our 
assumption that $a_1+a_2\in A_a$, 
we have $b+a_1+a_2\in A_b$. 
By 2 (applied to $b-a_2$), we have $b+a_1-a_2\in A_b$;
similarly, 
$$
b+2a_1, b+2a_1-a_2, b+2a_1+a_2\in A_b.
$$ 
By 3, applied to $b+a_1-a_2$,  
we have $b+3a_1\in A_a$.

\vskip 0,5cm
\includegraphics{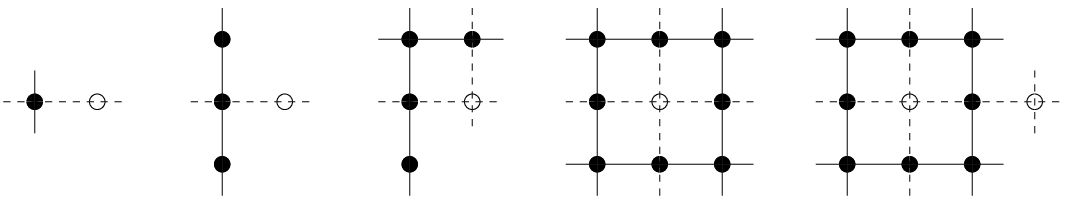}
\vskip 0,5cm

\noindent
Clearly, a pattern is emerging: we can rename $b+2a_1$ to 
$b$ and repeat the argument. 
Further, notice the symmetry with respect to $a_1$ and $a_2$, 
as well as the symmetry with respect to $\pm 1$. 
In conclusion, we have:

For any $a_3\in A_a$ 
such that $\langle a_1,a_2,a_3\rangle = A$ we have
$$
n_1a_1+n_2a_2 + a_3\in A_a
$$
iff both $n_1$ and $n_2$ are divisible by 2.  

\

5. We claim that for any primitive 
$x\in C$ and $c\in C$ we have $f(x+4c)=f(x)$.
Indeed, consider the lattice 
$$
E:=\langle x,a_3+2c\rangle.
$$
In $E$, we have two sublattices of index 2:
$$
\begin{array}{ccc}
E' & := & \langle a_3+2c,a_3+2c+2x\rangle\\
E''& := & \langle a_3+2c+x,a_3+2c-x\rangle.
\end{array}
$$
The generic values of $f$ on these sublattices
are different (by 4). It follows that one of them is equal $E^1_f$. 
By Corollary~\ref{coro:rk2-easy}, we have
$f(x+2y)=f(x)$ for every $y\in E$, in particular
$$
f(x+2a_3+4c)=f(x).
$$
Now consider the lattice 
$$
G:=\langle x', a_3\rangle,
$$
where $x'=x+4c$, and the sublattices 
$$
\begin{array}{ccc}
G' & := & \langle a_3, a_3+2x'\rangle\\
G''& := & \langle a_3+x', a_3-x'\rangle.
\end{array}
$$
Both have index 2 in $G$ and have different generic values. 
If follows that $|G/G^1_f|=2$. 
In particular, 
$$
f(x'+2a_3)=f(x+4c+2a_3)=f(x+4c).
$$
Combining with the result for $E$ we get our claim. 
It follows that for every sublattice $C\subset A$ of rank 2
$p(C)=2$ and, moreover, that $|C/C^1_f|$ is equal to 2 or 4. 
If $|C/C^1_f|=2$ for every subgroup $C$ of rank 2 
we get a contradiction to our assumption (this leads to Case 1). 

\

6. We can assume that $|C/C^1_f|=4$. We claim that $f$
is as in Case 2. 
First of all, 
$$
f(n_1a_1+n_2a_2+m_3a_3)=f(n_1a_1+n_2a_2+a_3)
$$
for all odd $m_3$ (equal to $f(a_3)$ iff $n_1=n_2=0$ modulo 2, by 4). 
Next, $f(2m_3a_3+c)=f(c)$ for all $m_3\in \Z$ and  all primitive  
$c\in C$. 
Since 
$$
f_C(x+4c)=f(x)
$$ 
for all primitive $x$ and all $c\in C$ we conclude
that either $f_C$ is constant, or induced from $C/2C$ or
as in Case 2. The first two possibilities contradict our 
assumptions on $C$.  
\end{proof}

\section{Reductions}
\label{sect:reductions}

\subsection{Reduction of $S$}
\label{sect:red-S}

\begin{lemm}\label{lemm:4coeff}
Let $f\in {\mathcal F}(A,S)$.
Assume that for all $h\,:\, S\ra \Z/4$
the function $h\circ f\in \AF(A,\Z/4)$.
Then $f\in \AF(A,S)$. 
\end{lemm}

\begin{proof}
The invariance 
is obvious (if $f(na)\neq f(a)$ for some 
$n,a$ then define $h$ so that $h\circ f(na)\neq h\circ f(a)$, 
leading to a contradiction).  
Assume that there exist elements $a,b\in A$ such that
$f(a), f(b)$ and $f(a+b)$ are pairwise distinct. 
Define $h$ such that 
$$
\begin{array}{lcc}
h\circ f(a)   & = & 0, \\
h\circ f(b)  & = & 1, \\
h\circ f(a+b)& = & 2.
\end{array} 
$$
Then, by Lemma~\ref{lemm:part-ord},  
$h\circ f\not\in \AF(A,\Z/4)$, 
contradiction. 
We see that for any $a,b\in A$ with 
$f(a)\neq f(b)$ either $f(a+b)=f(a)$ or $f(a+b)=f(b)$.
This defines a partial relation  $\tilde{>}$ on $A$
(as in Remark~\ref{rem:partial}).

We need to check that $\tilde{>}$ can be extended
to an order on $A$. 
Let $a,a'\in A$ be such that $f(a)=f(a')$. 
If there is no $b\in A$ such that $f(b)\neq f(a)$
then $f$ is constant and thus $\in \AF(A,S)$. 
If for all such $b\in A$ we have $a \tilde{>} b $ 
and $a'\,\tilde{>}\, b$ then 
$a=_f a'$. Otherwise, 
$b$ is a separator and we can assume
that $a\,\tilde{>}\, b\,\tilde{>}\, a'$. 
Assume that for some other separator $b'$
we have $a'\,\tilde{>}\, b'\,\tilde{>}\, a$. Let 
$$
\begin{array}{ccc}
h\circ f(a) & = & 0, \\
h\circ f(b) & = & 1 
\end{array}
$$ 
and put (if $f(b)\neq f(b')$) 
$$
h\circ f(b')=2.
$$ 
By assumption, 
$h\circ f\in \AF(A,\Z/4)$, contradiction (we use that
either $f(a)=f(a+b)$ or $f(b)=f(a+b)$, etc). 
Thus we have a correctly defined relation $>$ on $A$. 

Now we check the transitivity of $>$. 
Assume that we have elements $a,b,c\in A$ such that
$a>_f b> c$. Assume that $c\ge a$. If 
the values of $f$ on $a,b,c$ are pairwise distinct, 
put  
$$
\begin{array}{ccc}
h\circ f(a) & = & 0, \\
h\circ f(b) & = & 1, \\
h\circ f(c) & = & 2. 
\end{array}
$$
Since $h\circ f\in \AF(A,S)$ we get a contradiction. 
If $f(a)=f(b)$, let $a'$ be their separator; if 
$f(b)=f(c)$ let $b'$ be their separator and if $f(c)=f(a)$, 
let $c'$ be their separator: $c\ge c' \ge a$.    
Then there is a map $h\,:\, S\ra \Z/4$ such that
$h\circ f\not\in \AF(A,\Z/4)$, contradiction. 

Finally, we need to check that if $a> b$ and $a> b$, then
$a> b+c$. Again, we can introduce separators, if necessary, 
and proceed as above. 

To conclude we apply Lemma~\ref{lemm:part}.
\end{proof}

\begin{lemm}\label{lemm:z22}
Let $A$ be a finitely generated group,
$S$ a finite set and $f\in \cF(A,S)$. 
Assume that for all $h\,:\, S\ra \Z/2$
one has $h\circ f\in \AF(A,\Z/2)$.
Then $f\in \AF(A,S)$. 
\end{lemm}

\begin{proof}
As above, the invariance of $f$ is obvious. 
Following the proof of Lemma~\ref{lemm:4coeff}, 
observe that for all $a,b\in A$ with $f(a)\neq f(b)$ either
$f(a+b)=f(a)$ or $f(a+b)=f(b)$. Thus we have  
a partial relation $\tilde{>}$ 
on these pairs as in \ref{rem:partial}.  

Let $h\,:\, S\ra \Z/2$ be a nonconstant map and  
$$
S(h):=\{ s\in S\,|\, 
\exists\, a\in A\setminus A^1_{h\circ f}\,\,\,{\rm with}\,\,\,
f(a)=s\}.
$$
Let $h_0$ be a map such that 
$$
|S(h_0)|=\min_h(|S(h)|).
$$ 
We can assume that $S=\{1,...,n\}$ and that 
$S(h_0)=\{1,...,k_0\}$.

Assume that $1< k_0<n$.   
Let $a_1,...,a_{k_0}$ be  some elements
in $A\setminus A^1_{h_0\circ f}$ with $f(a_j)=j$. 
Then, for all $j\in S(h_0)$ and all $i\not\in S(h_0)$ we have
$a_j> x_i$ for all $x_i\in A$ with $f(x_i)=i$. 

Let $h'$ be the map sending each element in 
$\{n, 2,...,k_0\}$ to $0$ and each element
in $\{1,k_0+1,...,n-1\}$ to $1$.
One of the values is generic for $h'\circ f$.
Assume that $0$ is the generic value for $h'\circ f$. 
Then $a_n\notin A\setminus A^1_{h'\circ f}$
(indeed, if $a_n$ were generic then $a_n> x_1$ for
all $x_1$ with $f(x_1)=1$, contradiction to the previous). 
But then $|S(h')|\le k_0-1$, contradiction to 
the minimality of $k_0$. 

Assume that $1$ is the value of 
a generic element for $h'\circ f$.
Similarly,  the elements $a$ with $f(a)\in \{ k_0+1,...,n-1\}$
cannot be generic for $h'\circ f$. It follows that 
generic elements for $h'\circ f$ are mapped to $1\in S$. 
Contradiction to the assumption that $1<k_0$.

If $k_0=1$ then 
$$
A^1:=\{ a\in A\,|\, f(a)\neq 1\}
$$ 
is a proper subgroup and 
$f$ is constant on $A\setminus A^1$. 
Applying the same argument to $A^1$ we obtain a filtration 
$(A^n)$ such that $f$ is constant on 
$\overline{A}^n$ for all $n\in \N$. 
\end{proof}

\subsection{Reduction of the rank}
\label{sect:red-rank}

\begin{lemm}\label{lemm:restr3}
Let $A$ be an abelian group and 
$f\in  \cF(A,S)$. If for all subgroups $B\subset A$ 
with $\rk(B)\le 3$ the restriction $f_B\in \AF(B,S)$
then  $f\in \AF(A,S)$. 
\end{lemm}

\begin{proof}
By definition, it suffices to consider 
finitely generated $A$. 
Invariance of $f$ is clear. 
By Lemmas~\ref{lemm:4coeff} and ~\ref{lemm:z22}, 
it suffices to assume $S=\Z/2$. 
As in Remark~\ref{rem:partial}
and in the proof of Lemma~\ref{lemm:4coeff},
we can define a partial relation $\tilde{>}$ 
on $A$, which by assumption and by Lemma~\ref{lemm:part-ord} 
extends to an order on subgroups of $\rk \le 3$
(see Lemma~\ref{lemm:part-ord}). 
We will denote the induced order
on subgroups $C\subset A$ by $>_{C}$.
We need to show
that this order extends compatibly to $A$.  
Notice that for $C\subset D\subset A$
the order $>_{D}$ is stronger than the order 
$>_{C}$.

We have a decomposition of $A=A_a\cup A_b$ (preimages of
$0,1$).
As in the proof of Lemma~\ref{lemm:agf}, 
we will use the letter $a$ (resp. $b$) for elements in $A_a$ 
(resp. $A_b$). 
The proof is subdivided into a sequence of Sublemmas.

\begin{sublemm}(Correctness)
\label{sublemm:1}
There are no $a,a'$ and $b,b'$ such that 
$$
a\, \tilde{>}\, b\, \tilde{>}\, a'\, \tilde{>}\, b'\, \tilde{>}\, a.
$$
\end{sublemm}

\begin{proof}
Introducing the subgroups
$$
\begin{array}{cc}
C:= & \langle a,b,a'\rangle \\
D:= & \langle a',b',a\rangle \\
M:= & \langle b,b',a+a'\rangle \\
N:= & \langle b',a'+b,a\rangle 
\end{array}
$$
we obtain 
$$
a >_{C} b + a' >_{C} a' \,\,\, 
{\rm and}\,\,\, a' >_{D} b' + a >_{D} a.
$$ 
It follows that 
\begin{equation}\label{eqn:4}
f(a+a')=f(a+a'+b)=f(a),
\end{equation}
\begin{equation}\label{eqn:5}
f(a'+a+b')=f(a'),
\end{equation} 
\begin{equation}\label{eqn:7}
f(b+a')=f(b)
\end{equation}
By Equations~\ref{eqn:4} and \ref{eqn:5},
neither $b$ nor $b'$ can be generic in $M$. 
Thus $a+a'>_{M}  b,b'$ and 
\begin{equation}\label{eqn:6}
f(a+a'+b+b')=f(a+a')=f(a).
\end{equation} 
On the other hand, in $N$, the element $a$  
is not generic: $f(a+b')=f(b')$.
Since $f(b+a')=f(b)$ and  
$f(b + a' + a)=	f(a)$ (by ~\ref{eqn:4})
and since $a$ is not generic, the element $b+a'$ cannot
be generic. It follows that $b'$ is generic in $N$
and 
$$
f(b'+ (b+a')+a)=f(b'),
$$ 
contradiction to \ref{eqn:6}.   
\end{proof}

The sublemma implies that
we can extend $\tilde{>}$ to a relation  $>$ 
on the whole $A$.

\begin{sublemm}(Transitivity)
\label{sublemm:3}
If $x > y > z$ then $x > z$. 
\end{sublemm}

\begin{proof}
We have to consider 4 cases:

Case 1. $a > b > a'$. Transitivity follows from the definition. 

Case 2. $a > b > b'$.
Let $a'$ be the separator.  Assume $f(b'+a)=f(b)$. 
Then $b' > a$ and we have a 
contradiction to Sublemma~\ref{sublemm:1}. 
Thus $f(b'+a)=f(a)$ and $a > b'$. 

Case 3. $a > a' > b'$. Let $b$ be the separator. 
Again, if $b' > a$ we get a 
contradiction to Sublemma~\ref{sublemm:1}.

Case 4. $a > a' > a''$. Denote by $b$ the separator
between $a$ and $a'$.  
We have $a > b > a' > a''$. Applying case 2 (with 
$a$'s and $b$'s interchanged), we get
$b > a''$. Thus $a > a''$ (by the definition). 
\end{proof}

\begin{sublemm}(Additivity)
If $x,y,z\in A$  and $x > y$ and $x > z$ then $x > y+z$.
\end{sublemm}

\begin{proof}
We are looking at the following cases:

Case 1. $a > b,b'$.  
Neither $b$ nor $b'$ can 
be generic in the subgroup $\langle a,b,b'\rangle$. 
Thus $a$ is generic and the claim follows.

Case 2. $a > b$ and $a > b' > a'$. 
Then 
$$
f(b'\pm a')=f(b') \,\,\, {\rm and}\,\,\,
f(a+b'+a')=f(a).
$$ 

Case 2.1. $f(b+a')=f(b)$. Consider the subgroup 
$\langle a,b+a',b\rangle$. 
The element $b$ cannot be generic 
since $f(a+b)=f(a)$. It follows that $b+a'$ cannot be generic
since 
$$
f(b+a' - b )=f(a')\neq f(b+a').
$$ 
Thus $a$ is the generic element and $a > b+a'$. 

Case 2.2. $f(b+b')=f(b)$. Apply Case 1:
$a > (b+b') -(b'-a')$.  

Case 2.3. $f(b+a')=f(a)$, $f(b+b')=f(a)$. Consider
the subgroup $\langle b,b',a'\rangle$. 
The element $a'$ cannot be generic since $f(a'+b')=f(b')$. 
The element $b$ cannot be generic since $f(b+a')=f(a)$.
Finally, $b'$ cannot be generic since $f(b+b')=f(a)$. 
Contradiction.

Case 3. $a > b'> a'$ and $a > b''>_f a''$. 

Case 3.1. $f(b'+a'')=f(b')$.
Consider $\langle a,b',a''\rangle$:
$b'$ is nongeneric, therefore, $a''$ is also nongeneric - 
it follows that $a$ is generic and that 
$a > b'+a''$. Consider $\langle a,b',a'\rangle$: 
again $a$ is generic and $b',a'$ are nongeneric, thus
$a > (a'-b')$. Now we can apply case 1 or 2, depending
on the value of $f(a'-b')$. 

Case 3.2. By symmetry, we can assume that
both 
$$
f(b'+a'')=f(a'+b'')=f(a).
$$ 
Combining with the assumption of Case 3 we obtain
$$
a'' > b' > a' > b'' > a'',
$$
contradiction to Sublemma~\ref{sublemm:1}.
\end{proof}

\

Now we apply Lemma~\ref{lemm:part} to 
conclude the proof of Lemma~\ref{lemm:restr3}.
\end{proof}

\subsection{The exceptional case}
\label{sect:c3}

\begin{lemm}
\label{lemm:bad3}
Let $A=\Z\oplus \Z\oplus \Z$ and $f\in \cF(A,S)$
be a function such that 

$\bullet$ for every rank 2 sublattice $C\subset A$
we have $f_C\in \AF(C,S)$;

$\bullet$ $f\notin \AF(A,S)$;

$\bullet$ $f$ does not have a special basis.

Then $f$ takes exactly two values on $A$. Moreover, 
$f$ is either as in Example~\ref{exam:exam} or as
in Example~\ref{exam:exam2}.
\end{lemm}

\begin{proof}
We can assume that $f\,:\, A\ra S$ is surjective.
By Lemma~\ref{lemm:z22}, there exists an 
$$
h\,:\, S\ra \Z/2=\{ \circ, \bullet \}
$$ such that
$h\circ f\notin \AF(A,\Z/2)$.  
By Proposition~\ref{prop:special}, 
$h\circ f$ is either of the first or the second type. 

Let us consider the first case.  
By Example~\ref{exam:exam}, $h\circ f$ is given by

\vskip 0,5cm
\centerline{\includegraphics{bad2.eps}}
\vskip 0,5cm

Let $L$ be a line (=rank 2 lattice in $C$) 
reducing to the line through $(0:0:1)$ and 
$(0:1:1)$ modulo 2 and 
$P\in L$ be a point in  
$$
(h\circ f)^{-1}((0:0:1)). 
$$ 
By Corollary~\ref{coro:rk2-easy}, 
an AF-function takes only two values on a lattice of rank 2.
Since $f_L\in \AF(L,S)$ the value $f(P)$ is generic for $L$. 
Thus for every point $Q$ on $L$ in 
$$
(h\circ f)^{-1}((0:1:0))\,\,\,{\rm and}\,\,\,(h\circ f)^{-1}(0:0:1)     
$$
we have 
$$
f(P)=f(Q).
$$
Take any line $L'$
which modulo 2 passes through $(0:0:1)$ and $(1:1:0)$.
The value of $f$ on the point of intersection $L\cap L'$ is 
$f(P)$.  Since $f(P')$ is the generic value for $L'$
for every point $P'$ on $L'$ in the preimage of 
$h\circ f((0:0:1))$ we have $f(P')=f(P)$. Moreover,  
$f(P')=f(Q')$ for every $Q'$ in the preimage $h\circ f((1:1:0))$ in $L'$.
Therefore, for {\em any}  line $L$ (resp. $L''$) 
such that the reduction modulo 2 passes through 
$(0:0:1)$ and $(0:1:1)$ (resp. $(0:0:1)$ and $(1:1:0)$) 
the generic value is $f(P)$. In particular, for every point $R$ in 
the preimage of $(1:1:0)$ we have $f(R)=f(P)$.

Now consider the preimages of the points
$$
(1:0:1),\,\, (1:1:1),\,\, (0:1:1).
$$
Every one of those is generic for some rank 2 sublattice in $A$. 
Since these lattices intersect, we can apply the same
reasoning as above. It follows that $f$ can take only two values
on $A$ and, moreover, that $f$ is induced from $A/2A$.

The second case is treated similarly.  
\end{proof}

\subsection{Checking the AF-property}
\label{sect:checking}

We summarize the discussion of the previous sections:

\begin{prop}\label{prop:redu-p}
Let $A\in \cA_{q}$, $S$ a set and $f\in \cF(A,S)$.
Assume that $q>2$ and that 
for all subgroups $C\subset A$ of rank $\le 2$ one
has $f_C\in \AF(C,S)$. Then $f\in \AF(A,S)$.
\end{prop}

\begin{proof}
By Lemma~\ref{lemm:restr3} it suffices to consider the case
when $\rk(A)=3$. 
Assume that $f\notin \AF(A,S)$. 
By Lemma~\ref{lemm:z22}, there exists a map $h\,:\, S\ra \Z/2$ 
such that $h\circ f\notin \AF(A,\Z/2)$. By Lemma~\ref{lemm:red2-p},
there exists a subgroup $C\subset A$ of rank 2 such that
$h\circ f_C\notin \AF(C,\Z/2)$. In particular, $f_C\notin \AF(C,S)$.  
\end{proof}

\begin{prop}
\label{prop:import}
Let $k$ be a field of $\char(k)=0$ and $K/k$ an extension.  
Let $S$ be a ring such that $2s\neq 0$ 
for all $s\in S$. Assume that  $f\in \LF(K,S)$ and that
for all $\Z$-submodules $C\subset K$ 
of rank $\le 2$ 
one has $f_C\in \AF(C,S)$. Then 
$f\in \AF(K,S)$. 
\end{prop}

\begin{proof}
Assume otherwise. Then, by Lemma~\ref{lemm:restr3},
there exists a submodule $C\subset K$ 
of rank 3 such that $f_C\notin\AF(C,S)$. 
Moreover, there exists a map $h\,:\, S\ra \Z/2$ such that
$h\circ f_C\notin \AF(C,\Z/2)$. By Proposition~\ref{prop:red2-0},
we can assume that $C$ does not have a special basis
(by restricting to a proper subgroup). Then, 
by Proposition~\ref{prop:special}, $h\circ f_C$ 
has the form described in 
Example~\ref{exam:exam} or \ref{exam:exam2}.
In both cases,  
$f_C$ itself takes exactly two values. 

Up to addition of a constant (and
shifting $C$ using the logarithmic property, if necessary), 
we can assume that $f$ takes the values 0 and $s$
(for some $s\in S\setminus 0)$ and that 
$$
f(1)=0.
$$
Moreover, (in both cases!) 
there exist elements $x,y\in C\subset K$
such that 
$$
f(x)=f(y)=f(x+y)=f(x+y+1)=s,
$$
$$
f(x+1)=f(y+1)=0.
$$

\vskip 0,5cm
\centerline{\includegraphics{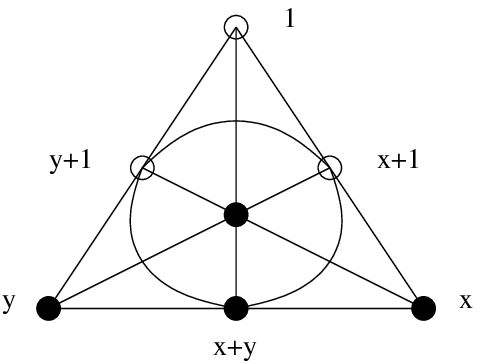}}
\vskip 0,5cm

\noindent
Consider the sublattices in $K$
$$
\begin{array}{ccc}
D & := & \langle xy,y,1\rangle;\\
E & :=& \langle xy,1,x+y\rangle.
\end{array}
$$
Using the logarithmic property we find that
$$
f(xy)=2s,\,\, f(y)=s, \,\, f(1)=0.
$$
By Lemma~\ref{lemm:bad3}, $f_D$ is an AF-function. 
Using the transitivity of the induced order on $D$ 
we see that 
$$
1 >_D y >_D xy.
$$ 
In particular, since $f(1)\neq f(xy)$, for any 
subgroup of $K$ containing  $1$ and $xy$ we have $1>xy$. 
Similarly, on $E$ the function $f_E$ also takes 3 values. 
Therefore (by Lemma~\ref{lemm:bad3}), 
$f_E\in \AF(E,S)$ and the induced order gives
$$
x+y>_E 1>_E xy
$$ 
(by transitivity of the order relation on $E$). 
It follows that 
$$
f(x+y+xy+1)=f(x+y)=s.
$$
On the other hand, (using the logarithmic property), 
$$
f(x+y+xy+1)=f((x+1)(y+1))=f(x+1)+f(y+1)=0.
$$
Contradiction.
\end{proof}

\section{AF-functions and geometry}  
\label{sect:afgeom}

\begin{assu}\label{assu:F}
Throughout $\cR$ is $\Q$, $\Z_p$
or $\Z/p$.
\end{assu}

\subsection{Affine geometry}
\label{sect:affine}

Let $k$ be a field  
and $V$ a (possibly infinite dimensional)
vector space over $k$, 
with an embedding of $k$ as $k\cdot 1$. 
For every pair 
$$
f_1,f_2\in \cF(\P(V),\cR)
$$ 
we have a map
$$
\begin{array}{cccc}
\varphi=\varphi_{f_1,f_2}: & \P(V)& \ra & \A^2(\cR)\\
                           &  v   & \mapsto & (f_1(v),f_2(v)).
\end{array}
$$

\begin{rem}\label{rem:cp}
If $f_1,f_2$ form a c-pair (see Definition~\ref{defn:c-pair}) 
then the image of every line in $\P(V)$ under $\varphi_{f_1,f_2}$ 
is contained in a line in $\A^2(\cR)$. 
\end{rem}

\begin{prop}
\label{prop:line-point}
If $f_1,f_2\in {\cF}(\P(V),\cR)$ 
form a c-pair then for every 
3-dimensional $k$-vector space $V'\subset V$
there exists an affine line ${\mathbb L}_{V'}\subset \A^2$
and a point $d_{V'}\in \A^2(\cR)$
such that 
$$
\varphi(\P(V'))\subset d_{V'}\cup {\mathbb L}_{V'}(\cR).
$$
\end{prop}

\begin{proof}    
Assume that $\varphi(\P(V'))$ contains 
4 distinct points $\overline{p}_1,...,\overline{p}_4$. 
Denote by $\overline{\ell}_{ij}$ the line through
$\ovl{p}_i$ and $\ovl{p}_j$.
Then the lines intersect in $\A^2(\cR)$ 
and the intersection point of 
these lines is contained in $\varphi(\P(V'))$.
Indeed, denote by 
$\ell_{ij}\subset \P(V')$ 
the line passing through
a pair of points $p_i, p_j$ contained
in $\varphi^{-1}(\ovl{p}_i)$ (resp. $\varphi^{-1}(\ovl{p}_j)$) 
and assume, for example, that 
$\ovl{\ell}_{12}\cap \ovl{\ell}_{34}=\emptyset$.   
The lines $\ell_{12}$ and $\ell_{34}$ 
intersect (or coincide) in $\P(V')$. 
The point of intersection is 
contained in the image. 
Thus the lines $\ovl{\ell}_{12},\ovl{\ell}_{34}$
intersect and the image of $\P(V)$ contains 
the intersection point.

Now we can assume that $\varphi(\P(V'))$ contains at 
least 5 points $\ovl{p}_j$, such 
that $\ovl{p}_1,\ovl{p}_2,\ovl{p}_3\in 
\overline{\ell}\in \A^2_{\cR}$
and $\ovl{p}_4,\ovl{p}_5\not\in\overline{\ell}$ and
such that the line through  $\ovl{p}_4,\ovl{p}_5$
intersects $\ovl{\ell}$ in $\ovl{p}_3$.

\vskip 0,5cm
\hskip 3cm\includegraphics{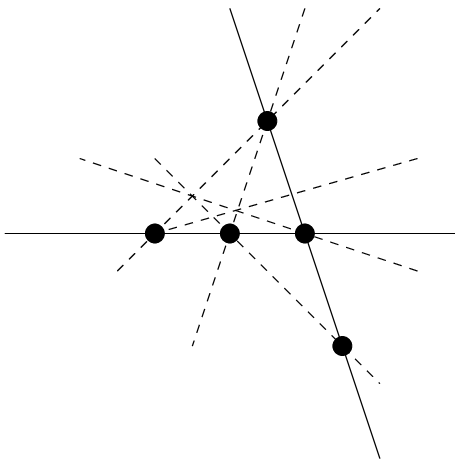}
\vskip 0,5cm

The goal is to show that drawing lines in $\A^2$  
through the points of intersections of the existing lines
one can produce
2 new pairs of points in $\A^2(\cR)$ such that 
the corresponding lines are parallel, leading to a 
contradiction.

First we assume that $\cR$ is $\Q$ or $\Z/p$.
Compactify $\A^2$ to $\P^2$ by adding a line $\ell_{\infty}$. 
After a projective transformation of $\P^2$ 
(with coefficients in $\cR$)
we can assume that the points are given by
$$
(1:0:1), \, (0:0:1),\, (0:1:1),\, (1:0:0),\, (0:1:0).
$$  
We use standard affine 
coordinates $(z_1,z_2)$ on $\A^2$ 
and corresponding projective coordinates $(z_1:z_2)$ 
on $\ell_{\infty}$.  The set $\varphi(\P(V'))$ has the
following properties:

$\bullet$ if $(z_1,z_2)\in \varphi(\P(V'))$ 
then $(z_1:z_2)\in \ell_\infty$
is also contained in $\varphi(\P(V'))$.

$\bullet$ if $(z_1,z_2)\in \varphi(\P(V'))$ 
then $(z_1,0),\, (0,z_2)$ are also
in $\varphi(\P(V'))$.

$\bullet$ if $(z_1,\,z_2) \in \varphi(\P(V'))$ 
then $(z_1:z_2:0)\in \varphi(\P(V'))$ and
$(z_1/z_2,0)\in \varphi(\P(V'))$.

$\bullet$ if $(z_1,0),\, (z_2,0)\in \varphi(\P(V'))$ then
$(z_1+z_2,0)\in \varphi(\P(V'))$.

$\bullet$ if $(z,0)\in \varphi(\P(V'))$ 
then $(0,z) \in \varphi(\P(V'))$.

To check the listed properties it suffices to compute
the coordinates of the points of intersection of
appropriate lines. For example, the first property
follows from the fact that
$(z_1:z_2)\in \ell_\infty$ is the
intersection of $\ell_\infty$ with the line 
through $(0,0)$ and $(z_1,z_2)\in \A^2$.
For the second, observe that $(z_1:0:1)$ can be obtained as 
the intersection of the line through
$(0:0:1) $ and $(1:0:1)$ with the line 
through $(z_1:z_2:1)$ and $(0:1:0)$, etc.

In particular, $\varphi(\P(V'))$  contains a subset 
of points $(z_1,z_2)$, where $z_1,z_2$ are generated 
from coordinates of points in 
$\varphi(\P(V'))\cap (\P^2\setminus \ell_\infty)$
by the above procedures. 
For $\cR=\Q$ or $\cR=\Z/p$, the set $\varphi(\P(V'))$ 
contains $\A^2(\cR)\cup\ell_{\infty}(\cR)=\P^2(\cR)$.  
In particular, one can find two lines in $\P(V')$ 
such that their images 
don't intersect in $\A^2$, contradiction. 

Now we show how to extend this argument to $\cR=\Z_p$. 
As before, we assume that 
$\varphi(\P(V'))$ contains 5 points as in the picture above.
One can choose a coordinate system  such that
$\varphi(\P(V'))$ contains the points $(1,0), (0,0), (0,1), (z_1,0)$
and $(0,z_2)$ with $z_1,z_2\in \Q_p$. 
Then it also contains some point $(z_1',z_2')$ with {\em nonzero} 
coordinates
$z_1',z_2'\in \Q_p$.
The (projective) transformation $\cT$ moving 
$$
(1,0), (0,0), (0,1)\,\, \mapsto
(\infty,0), (0,0), (0,\infty)
$$
is given by
$$
(w_1,w_2)=(\frac{z_1}{1-(z_1+z_2)}, \frac{z_2}{1-(z_1+z_2)})
$$
and its inverse $\cT^{-1}$ by 
$$
(z_1,z_2)=(\frac{w_1}{1+w_1+w_2}, \frac{w_2}{1+w_1+w_2}).
$$
Apply the reasoning of the first part to 
$(w_1',w_2'):=\cT((z_1',z_2'))$.
First we find that $\varphi(\P(V'))$ contains the points
$(w_1',0), (0,w_2')$, then that it contains all points of the form
$(w_1'/2^{m_1},w_2'/2^{m_2})$ for some $m_1,m_2\in \N$, then 
that it contains all points $(r_1w_1'/2^{m_1},r_2w_2'/2^{m_2})$
with $r_1,r_2\in \N$, and, finally,
that it contains all points  
with coordinates $(e_1w_1',e_2w_2')$ 
for arbitrary $e_1,e_2\in \Q$. 

To arrive at a contradiction it suffices to produce a pair of rational
numbers $e_1,e_2$  
such that 
$$
\cT^{-1}((e_1w_1',e_2w_2'))\not\in \Z_p\oplus \Z_p.
$$
Clearly, (for any $w_1',w_2'\in \Q_p$)  we can find $e_1,e_2\in \Q$ 
such that
$$
\frac{e_1w_1'}{1+e_1w_1'+e_2w_2'}\not\in \Z_p.
$$
This concludes the proof.
\end{proof} 

\begin{rem}
Proposition~\ref{prop:line-point} is wrong 
for  $\cR=\Q_p$. 
\end{rem}

\subsection{Projective geometry}
\label{sect:geometry}

\begin{prop}\label{prop:2-coeff}
Let $k$ be any field and
$V$ a 3-dimensional vector space over $k$. Assume that 
$f_1,f_2\in {\cF}(\P(V),\Z/2)$ are
such that 
$$
\varphi(\P(V))\subset \{ \ovl{p}_{12},\ovl{p}_{13},
\ovl{p}_{23}\}\subset \A^2(\Z/2),
$$ 
where $\ovl{p}_{12}=(0,0),
\ovl{p}_{13}=(1,0)$ and $\ovl{p}_{23}=(0,1)$. 
Assume further that for all 
2-dimensional vector spaces 
$V'\subset V$ the image $\varphi(\P(V'))$ is contained
in at most two of these points. 
Then at least one of the functions 
$f_1,f_2$ or $f_3=f_1+f_2$ is
an AF-function on every 2-dimensional subspace $V'\subset V$.
\end{prop}

\begin{nota}
Denote by $P_{ij}\subset \P(V)$ the preimage of $\ovl{p}_{ij}$.
Let $T_{i}$ be the 
set of lines $t\subset \P(V)$ 
such that $\varphi(t)\subset \{ \ovl{p}_{ij}, \ovl{p}_{ik}\}$. 
\end{nota}

The proof of  Proposition~\ref{prop:2-coeff}
is subdivided into a sequence of lemmas.

\begin{lemm}\label{lemm:picture}
Let $t_{i}, t_{i}'$ be two lines in $T_{i}$.
Every point $p_{jk}$ defines a projective isomorphism between
$t_{i}\cap P_{ij}$ and $t_i'\cap P_{ij}$. 
\end{lemm}

\begin{proof}

\vskip 0,5cm
\hskip 2cm\includegraphics{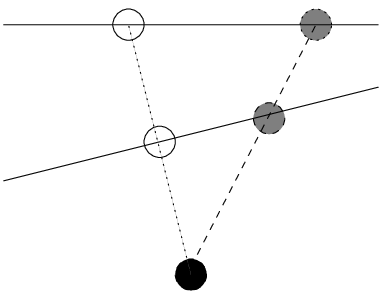}
\vskip 0,5cm

\end{proof}

\begin{lemm}\label{lemm:except}
If there exists a line $t_{i}\in T_{i}$
such that the number of points in 
$t_{i}\cap P_{ij}$ is $\le 1$ then 
either $f_1,f_2$ or $f_3$ is 
a GF-function on $V$. 
\end{lemm}

\begin{proof}
First of all, if one of the sets $P_{ij}$ is empty
then one of the functions $f_1,f_2,f_3$ is the constant
function, hence a flag function.

Assume that this is not the case and 
that there exists a line $t_i$ such that 
$t_{i}\cap P_{ij}=\emptyset$. 
By assumption, there exist points of 
all three types. We can draw a line $t_j$
through some points of type $P_{ij}$ and $P_{jk}$. 
The line $t_j$ intersects $t_i$ in a point, 
which must be a point of type $P_{ik}$. 
Thus the line $t_i$ 
contains points of all three types, contradiction.

Finally, assume that there exists a line $t_i$ such 
that $t_{i}\cap P_{ij}$ consists of exactly one point. 
There are two possibilities:
there are at least two lines of type $T_j$ 
or there is exactly one line of type $T_j$.
In the first case Lemma~\ref{lemm:picture}, shows that 
{\em all} lines of type $T_j$ contain exactly one point
of type $P_{ij}$. This means that there is only one point of type
$P_{ij}$ in $\P(V)$ (otherwise, we could draw a line through two
of those points; this line cannot 
by of type $T_i$ nor of type $T_j$). 
It follows that $f_k$ is an 
AF-function on $V$ 
(delta function). 

Assume now that there exists exactly one line of type $T_j$.
The complement to this line contains {\em only} points of type
$P_{ik}$. It follows that $f_j$ is an AF-function on $\P(V)$. 
\end{proof}

Thus we can assume that there are at least 3 
points of each type (which do not lie on a line) 
and that there are at least two lines of each type and that
for every line  $t_{i}\in T_{i}$
the set $t_{i}\cap P_{ij}$ has at least two elements.

\begin{defn}
\label{defn:rel}
We call the points $p_{ij},p_{ij}'\in t_i$ {\em related} if 
there exists a point $p_{ij}^0$ such that the line
joining $p_{ij}^0$ and $p_{ij}$ and the line
joining $p_{ij}^0$ and $p_{ij}'$ are both of type $T_j$. 
\end{defn}

\begin{lemm}\label{lem:14.4}
If there exists a line $t_i\in T_i$ containing 
two nonrelated points $p_{ij},p_{ij}'\in t_i$ then 
for every $t_j\in T_j$ passing through 
$p_{ij}$ or $p_{ij}'$ all points in $t_{i}\cap P_{ij}$ are related. 
\end{lemm}

\begin{proof}
Consider this line $t_i$ with 
two nonrelated points $p_{ij},p_{ij}'$. Let $t_j\in T_j$ be any line
passing through $p_{ij}$. Let $p_{ij}''$ be an arbitrary 
point in $t_{i}\cap P_{ij}$, distinct from $p_{ij}$. 
Since  $p_{ij},p_{ij}'$ are not related the line through
$p_{ij}''$ and $p_{ij}'$ has to be of type $T_i$. 
It follows that all points of type $P_{ij}$ on $t_j$
are related through $p_{ij}'$.  
\end{proof}

\begin{lemm}\label{lemm:auto}
Assume that $p_{ij},p_{ij}'\in t_i$ are related. 
For every point $p_{ik}\in t_i$ there exists a
projective automorphism $m_{ik}\,:\, t_i\ra t_i$ such that

$\bullet$\,\, $m_{ik}(t_i\cap P_{ij})=t_i\cap P_{ij}$;

$\bullet$\,\, the unique fixed point of $m_{ik}$ is $p_{ik}$;

$\bullet$\,\, $m_{ik}(p_{ij})=p_{ij}'$.
\end{lemm}

\begin{proof}
Consider the triangle 
spanned by $p_{ij},p_{ij}'\subset t_i$ and $p_{ij}^{0}$. 
Draw a line $t_i'$ through $p_{ik}\in t_i$ and $p_{ij}^0$. 
Pick a point $p_{jk}$ on the line $t_j$ through
$p_{ij}$ and $p_{ij}^0$ and draw a line $t_k$ through
$p_{jk}$ and $p_{ik}$. 
Denote by $p_{jk}'$ the point of intersection 
of $t_j$ with the line through $p_{ij}'$ and $p_{ij}^0$. 
Using the points $p_{jk}$ and $p_{jk}'$ as centers for 
the projective isomorphism between the line joining 
$t_i$ and the line $t_i'$, we obtain the projective 
automorphism $m_{ik}$.

\vskip 0,5cm
\hskip 3cm\includegraphics{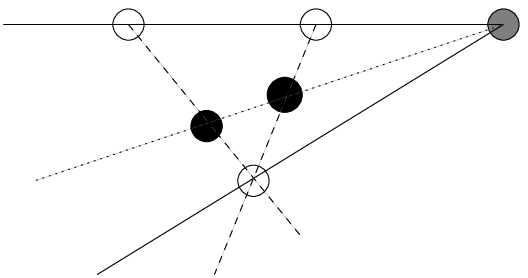}
\vskip 0,5cm

\end{proof}

\begin{lemm}\label{lemm:ff}
Assume that there exists a line $t_i$ such that
all points $p_{ij}\in t_i$ are related. Then all three
functions $f_1,f_2,f_2$ are AF-functions on all lines of 
type $T_i$. 
\end{lemm}

\begin{proof}
The function $f_i$ is constant on lines of type $T_i$, 
hence an AF-function. 
The function $f_j+f_k$ is constant on (the fixed line) $t_i$. 
Let us show that $f_j$ is
an AF-function on $t_i$. 

By Lemma~\ref{lemm:auto}, for any pair
of related points $p_{ij},p_{ij}'$ and any point $p_{ik}$  on 
a projective line of type $t_i$
there exists a projective automorphism (transvection)
with a single fixed point $p_{ik}$ moving $p_{ij}$ to $p_{ij}'$.
Introduce coordinates on $\P^1$ such that $p_{ik}=(0:1)$,
$p_{ij}=(1:0)$ and $p_{ij}'=(1:1)$. 
A unipotent lifting of the automorphism $m_{ij}$ 
to $\GL(V)=\GL_2(k)$ 
can be written in the form
$$
\left(\begin{array}{cc} 1 & 1\\ 0 & 1\end{array}\right),
$$
In this basis $f(e_1)=f(e_1+e_2)=1$ and $f(e_2)=0$. 
Consequently, $f$ satisfies the conditions 
of Lemma~\ref{lemm:agf}.
(Notice that the points $p_{ik}$ on $t_i$
need not be related. This leads to the nonsymmetric
shape of the functional equation~\ref{eqn:fe}.)
Now we apply Lemma~\ref{lemm:agf}.

To conclude that $f_1,f_2,f_3$ are AF-functions on
{\em all} lines of type $T_i$ we use
the projective isomorphism preserving 
both sets $t_{i}\cap P_{ij}$ and 
$t_{i}\cap P_{ik}$ (introduced in Lemma~\ref{lemm:picture}).
\end{proof}

\begin{proof}(of Proposition~\ref{prop:2-coeff})
We may assume that we are not in the situation of 
Lemma~\ref{lemm:except}. 
If for all $i=1,2,3$ and all lines $t_i$ of type $T_i$
all points in $t_{i}\cap P_{ij}$ or
all points in $t_{i}\cap P_{ik}$ are related then 
all three functions $f_1,f_2,f_3$ are 
AF-functions on $t_i$.  

Assume there is a $t_i$ and two points 
$p_{ij}, p_{ij}'\in t_{i}\cap P_{ij}$
which are not related. By Lemma~\ref{lem:14.4} and
Lemma~\ref{lemm:ff}, $f_1,f_2,f_3$ are AF-functions on 
all lines of type $T_j$. There are two cases: 
there exist two nonrelated points in $t_{i}\cap P_{ik}$ or not.
In the first case, $f_1,f_2,f_3$ are  AF-functions on 
all lines of type $T_k$. In the second case $f_1,f_2,f_3$ 
are AF-functions on lines of type $T_i$. 
If, for example, all three functions 
are AF-functions on lines of type
$T_j$ and $T_k$ then the function $f_i$ 
(being constant on lines of type $T_i$)
is an AF-function on {\em all} lines in $\P(V)$.
This concludes the proof.
\end{proof}

\subsection{Logarithmic functions}
\label{sect:logo}

We keep the assumptions of \ref{assu:F}.

\begin{prop}
\label{prop:lambda}
Let $k$ be a field of characteristic $\neq 2$
and $K/k$ an extension. Assume that  
$f_1,f_2\in \LF(K,\cR)$ 
form a c-pair (see \ref{defn:c-pair} for the definition). 
Assume that the linear space $\langle f_1,f_2\rangle_{\cR}$
does not contain a (nonzero) AF-function. 
Then there exists a 3-dimensional $V\subset K$ such 
that for {\em every} (nonzero) 
$f'\in \langle f_{1,V},f_{2,V}\rangle_{\cR}$
we have $f'\notin\AF(V,\cR)$.    
\end{prop}

\begin{proof}
For $\char(k)>2$ we use Proposition~\ref{prop:redu-p}
and for $\char(k)=0$ Proposition~\ref{prop:import}.
Since $f_1\notin\AF(K,\cR)$ 
there exists a 2-dimensional subspace  $V'\subset K$ 
such that
$f_{1,V'}\notin \AF(V',\cR)$. 
Since 
$$
\rk \, \langle f_{1,V'},f_{2,V'},1\rangle \le 2
$$  
we have
$f_{2,V'}-\mu_1 f_{1,V'}=\mu_2$ (for some $\mu_1,\mu_2\in \cR$).  
Since 
$$
f_2-\mu_1 f_1\not\in \AF(K,\cR),
$$
by Section~\ref{sect:checking}, 
there exists a 2-dimensional $W'$
such that 
$$
f_{2,W'}-\mu_1 f_{1,W'}\not\in \AF(W',\cR).
$$ 
Choose some $k$-basis in 
$V'=\langle x_1,x_2\rangle$ and $W'= \langle y_1,y_2\rangle$
(with $x_j,y_j\in K^*$). 
Let $V=\langle x_1,x_2, y_2y_1^{-1}x_1\rangle$. 
Then for every pair of $(\la_1,\la_2)\neq (0,0)$ 
$$
\lambda_1 f_1 +\lambda_2(f_{2}-\mu_1 f_{1})\not\in \AF(V,\cR).
$$
Indeed, for pairs with 
$\la_1\neq 0$ consider the restriction to $V'$.  
For pairs $(0,\la_2)$ with $\la_2\neq 0$
consider the restriction to (a shift of) $W'$ and 
use the invariance and the logarithmic property of $f$. 
\end{proof}

\begin{lemm}\label{lemm:h}
Let  $k$ be a field of characteristic
$\neq 2$, $K/k$ an extension 
and $V\subset K$ a 3-dimensional vector space over $k$. 
Consider a c-pair  $f_1,f_2\in \cF(\P(V),\cR)$. 
Assume that there are no nonzero
AF-functions $f\in \langle f_1,f_2\rangle_F$.
Then there exist nonconstant (and nonproportional) functions 
$f_1', f_2'\in  \langle f_1,f_2,1\rangle_{\cR}$
and a map $h'\,:\, \cR\ra \Z/2$
such that neither of the three functions
$$
\begin{array}{c}
h'\circ f_1',\\
h'\circ f_2',\\
h'\circ f_1'+h'\circ f_2'
\end{array}
$$ 
is an AF-function on $V$. 
\end{lemm}

\begin{proof}
By \ref{prop:line-point}, we know that
$$
\varphi_{f_1,f_2}(\P(V))\subset d_{V}\cup {\mathbb L}_{V}.
$$ 
After a linear change of coordinates (over $\cR$) 
we can assume that
$$
\varphi_{\tilde{f}_1,\tf_2}(\P(V))=(0,1)\cup \{ x-{\rm axis}\},
$$
where 
$$
\begin{array}{ccc}
\tf_1 & = & \la_1f_1+\la_2f_2+\la_3, \\
\tf_2 & = & \mu_1f_1+\mu_2f_2+\mu_3
\end{array}
$$
and

$\bullet$ $\tf_1(0)=\tf_2(0)=0$;

$\bullet$ $\tf_2$ takes only two values;

$\bullet$ $\tf_2(v)=0$ if $\tf_1(v)\neq 0$;

$\bullet$ $\tf_1(v)=0$ if $\tf_2(v)\neq 0$.

\noindent
Let $h$ be such that $h\circ \tf_1\notin \AF(V,Z/2)$. 
We can assume that $h(0)=0$. 
Let $v_1\in V$ be such that $h\circ \tf_1(v_1)\neq 0$. 
After rescaling (and a corresponding rescaling of $h$), 
we can assume that $\tf_1(v_1)=1$ and that $h(1)=1$. 
Since $\tf_2$ takes only the values $0,1$
we have 
$$
h\circ \tf_2=\tf_2 \notin \AF(V,Z/2).
$$

Notice that the last two properties imply that
\begin{equation}\label{eqn:h}
h\circ (\tf_1+\tf_2) = h\circ \tf_1+h\circ\tf_2.
\end{equation}

The lines in $\P(V)$ can be subdivided into 3 classes
according to their image under $\varphi_{\tf_1,\tf_2}$.

\vskip 0,5cm
\hskip 3cm\includegraphics{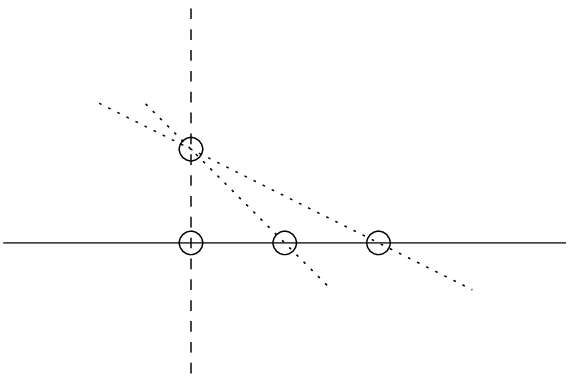}
\vskip 0,5cm

The function $\tf_1$ is obviously constant 
on lines in $\P(V)$ whose image is contained in the $y$-axis, 
$\tf_2$ is constant on lines mapping to the 
$x$-axis and 
$\mu \tf_1+ \tf_2$ (with $\mu\neq 0$) is constant
on the lines mapping to $\mu x+ y = 1$.  
We call the lines of the first type $T_1$-lines, the lines
of the second type $T_2$-lines and all other $T_3$-lines.

Notice that $\tf_3$ restricted to the $T_i$-lines coincides with $\tf_i$
for $i=1,2$.  
By assumption, $\tf_3\notin\AF(V,\cR)$, and by the results 
in Section~\ref{sect:checking}, 
there exists a line $\ell_3 \subset \P(V)$ such that $\tf_3\notin
\AF(\ell_3,\cR)$. If $\ell_3$ is of type $T_1$ or $T_2$ 
then $h$ as above solves our problem.
Thus we can assume that 
all three functions $\tf_1,\tf_2$ and $\tf_3$ are AF-functions on all lines
of type $T_1$ and $T_2$. 
Consider a line $\ell_3$ of type $T_3$ such that 
$\tf_3\notin \AF(\ell_3, \cR)$. 
Then $\varphi(\ell_3)\subset \mu x+ y = 1$, (with  $\mu\neq 1$). 
That is, on $\ell_3$ we have 
$$
\begin{array}{ccc}
\mu \tf_1+\tf_2 & = & 1\\
\tf_1+\tf_2     & = & \tf_3
\end{array}
$$ 
It follows that all 3 functions $\tf_1,\tf_2,\tf_3\notin\AF(\ell_3,\cR)$. 
Now we change coordinates again, making $\ell_3$ the new coordinate
axis (in addition to the $x$-axis).  
We put 
$$
\begin{array}{ccl}
f_1' & = & \mu \tf_1+\tf_2 - 1\\
f_2' & = & -\tf_2  \\
f_3' & = &  f_1'+f_2'.
\end{array}
$$ 
Now we can apply the argument above: find an $h'$ such that 
$$
h'\circ f_1'\notin \AF(V,\Z/2);
$$ 
(after rescaling, if necessary) 
one can assume that 
$$
h'\circ f_2'=f_2',
$$ 
(with $f_2'=\la \tf_2$ for some nonzero $\la$).
Since $(f_1'+f_2')$ restricted to 
$\ell_3=\{ f_1'=0\}$ is equal to $f_2'$, which is not AF on $\ell_3$,
(by our assumption that $\tf_2$ is not AF on $\ell_3$), 
we conclude that  
$$
h'\circ (f_1'+f_2')\notin \AF(\ell_3, \Z/2).
$$ 
By the same reasoning as above (see \ref{eqn:h}) we have
$$
h'\circ (f_1'+f_2')= h'\circ f_1'+ h\circ f_2'.
$$
This finishes the proof. 
\end{proof}

\subsection{Existence of AF-functions}
\label{sect:existence}

\begin{prop}
\label{prop:cp}
Let $K/k$ be an extension of fields.
Assume that $f_1,f_2\in \LF(K,\cR)$ form a c-pair.
Then $\langle f_1,f_2\rangle_{\cR}$ contains an AF-function. 
\end{prop}

\begin{proof}
We can assume that $f_1,f_2,1$ 
are linearly independent on $\P(K)$.
(Otherwise some linear combination  
is constant, hence an AF-function.)

Assume that $\langle f_1,f_2\rangle_{\cR}$ does not 
contain an AF-function.
By Proposition~\ref{prop:lambda} combined with  
Lemma~\ref{lemm:h}, there exist 
a 3-dimensional $V\subset K$, functions 
$f_1',f_2'\in \cF(\P(V),\cR)$
and a map $h\,:\, \cR\ra \Z/2$ such that
$$
h \circ f_1', h\circ f_2', 
h \circ f_3'= h\circ f_1'+h\circ f_2'\notin \AF(V,\Z/2).
$$
These functions satisfy the assumptions of
Proposition~\ref{prop:2-coeff}. 
We obtain a contradiction to its statement. 
\end{proof}

\section{Galois theory}

\subsection{Groups}
\label{sect:groups}

Let $G$ be a (topological) group with unit $0=0_G$ 
and $g,g'\in G$. 
Denote by  $[g,g']$ their commutator and by 
$$
G=G^{(0)}\supset G^{(1)} \supset ... 
$$
the lower central series: $G^{(i)}$ is the (closed) 
subgroup generated by 
$[g_i,g_0]$, where $g_i\in G^{(i)}$ and $g_0\in G^{(0)}$. 
Denote by $G^a=G/G^{(1)}$ the abelianization of $G$,
by $G^c=G/G^{(2)}$ the second quotient - it is a 
central extension of $G^a$ - and by  
$G^{1,2}=G^{(1)}/G^{(2)}$. 
Let
$$
\psi_c\,:\, G\ra G^c, \hskip 1cm 
\psi_a\,:\, G\ra G^a, \hskip 1cm 
\psi_{ca}\,:\, G^c\ra G^a
$$
be the quotient homomorphisms 
(we have $\psi_a=\psi_{ac}\circ\psi_c$).

\begin{lemm}
\label{sect:seq}
Let $G$ be a (profinite) group and 
$\psi\,:\, G\ra A$ a (continuous) surjective homomorphism 
onto a finite group $A$. Assume that 
$\al\in H^2(A,\Z/p^n)$ is a class such that 
$$
\psi^*(\al)=0\in H^2(G,\Z/p^n).
$$ 
Then 
there exists a (continuous) homomorphism 
$\tilde{\psi}\,:\, G^c\ra A$ such that
$\tilde{\psi}\circ \psi_c=\psi$
and 
$$
\tilde{\psi}^*(\al)=0\in H^2(G^c,\Z/p^n).
$$
\end{lemm}

\

We write $G(p^n)$  for the subgroup of 
$G$ generated by 
$g^{p^n}$ with $g\in G$ 
(for abelian groups and their central
extensions, the product of 
$p^n$-th powers is a $p^n$-th power).     
For any  profinite group $G$ we denote
by $G_p$ its maximal pro-$p$-quotient.
Notice that for any profinite group $G$ we have 
$(G^a)_p=(G_p)^a$ and $(G^c)_p=(G_p)^c$.

\subsection{Fields}
\label{sect:fields}

\begin{assu}\label{assu:fff}
Fix a prime $p$. We assume that $\char(k)\neq p$ and that
$k$ does not admit finite separable extensions of
degree divisible by $p$. 
\end{assu}

Let $K$ be a field over $k$.
It has a structure of a vector space
over $k$ and therefore, $K^*/k^*$ 
a structure of a projective
space over $k$, though infinite-dimensional. 
We continue to denote this space by $\P(K)$ and by
${\cF}(\P(K),\Z_p)$
the space of $\Z_p$-valued 
functions on $\P(K)$. 

In Section~\ref{sect:defn} we defined the set
$\LF(K,\Z_p)$. We now consider the topological space
$\LF^{\rm top}(K,\Z_p)$. As a set it coincides with
$\LF(K,\Z_p)$. The basis of the
topology is given by $U_{E^*,n}$, where
$E^*$ is a finitely generated subgroup of $K^*/k^*$ and 
$n\in \N$ - a function $f$ is in $U_{E^*,n}$
if $f$ is equal to $0$ modulo $p^n$ on $E^*$.

\begin{nota}
\label{nota:ff}
For any $k$-vector space $V\subset K$
(not necessarily closed under multiplication
in $K$)
and $f\in {\cF}(\P(K),\Z_p)$ we denote by $f_V$
the restriction of $f$ to
$\P(V)=(V\setminus 0)/k^*$ (sometimes 
we will denote by the same 
symbol the restriction of $f$ to $V$). 
For a finite 
set of functions $f_j\in {\mathcal F}(\P(K),\Z_p)$ 
we denote by $\langle f_1,...,f_n\rangle $ the 
$\Q_p$-vector space they span in 
${\cF}(\P(K),\Z_p)_{\Q_p}$.   
\end{nota}

\subsection{Galois groups}
\label{sect:galois-groups}

Let $K$ be a field as in Section~\ref{sect:fields}.   
Denote by $\Gal_K$ the 
Galois group of a separable closure of $K$.
It is a profinite compact topological group (we refer to 
\cite{serre} for basic facts concerning 
profinite groups). 
In general, the group $\Gal_K$ has a rather complicated 
structure. We will be interested in 
$$
\Gal_{K/k}:=\Ker(\Gal_K\ra \Gal_k),
$$
more precisely, in the pro-$p$-group 
$\G^c:=(\Gal^c_{K/k})_p$
and its abelianization $\G^a:=(\Gal^a_{K/k})_p$.
The group  $\G^a$ is a torsion free abelian pro-$p$-group
(by our assumptions on $k$).

\begin{lemm}\label{lemm:ga-functions}
One has a (noncanonical) isomorphism of topological groups
$$
\G^{a}\simeq \LF^{\rm top}(K,\Z_p).
$$
\end{lemm}

\begin{proof}
Since $K$ contains all $p$-power roots of $1$ 
we can choose an identification between
$\Z_p$ and $\Z_p(1)$. Then, by
Kummer theory, we have a nondegenerate pairing 
$$
\G^a/\G^a(p^n)\times K^*/(K^*)^{p^n}\ra \Z/p^n,
$$
given by 
$$
(\sigma,\kappa)\mapsto \sigma(\kappa)/\kappa
$$
for $\sigma\in \G^a/\G^a(p^n)$ and $\kappa\in K^*$.
We derive an
isomorphism (of topological groups)
$$
\G^a = \Hom(\hat{K}^*,\Z_p)
$$ 
(where $\hat{K^*}$ is the completion of $K^*$ with 
respect to subgroups of $p$-power index).
Moreover, every such
homomorphism is trivial on $k^*$, since $k^*$ 
does not admit finite extensions of degree divisible 
by $p$, by assumption.
Thus, in our case, the latter group is isomorphic to 
$\LF^{\rm top}(K,\Z_p)$, 
(since $K^*$ is dense in $\hat{K}^*$).

More explicitely, 
two elements of the Galois group $\G^a$ 
coincide if for all $n$ their actions 
on all cyclic $p^n$-degree extensions of 
$K$ coincide.   
Thus the resulting map 
$$
\G^a/\G^a(p^n) \ra \LF^{\rm top}(K,\Z/p^n)
$$
is a monomorphism.

Conversely, every element of
$\LF^{\rm top}(K,\Z_p)$
defines an element of $\G^a$.  
Any homomorphism 
$\chi\,:\, K^*/k^*\ra \Z_p$ defines
a compatible set of elements of 
$\G^a/\G^a(p^n)$ for all abelian
extensions of $K$ of degree $p^n$. 
Thus $\G^a$ and $\LF^{\rm top}(K,\Z_p)$ 
are isomorphic as topological
groups.

(We use the fact that the group $K^*/k^*$ 
has no torsion, by assumptions on $k$.  
Therefore, the map  
$$
\LF^{\rm top}(K,\Z/p^{n+1})\ra \LF^{\rm top}(K,\Z/p^n)
$$ 
corresponding to the projection 
$\Z/p^{n+1}\to \Z/p^n$ is surjective and 
we obtain an isomorphism between projective limits.)
\end{proof}

We  consider $K$ 
as a vector space over $k$, with
a canonical embedding of 
$k$ as a 1-dimensional subspace $k\cdot 1$.  
Every finite-dimensional subspace $V\subset K$
which contains $k$ (as a subspace) 
defines a subfield  $K_V$ of $K$
(generated by elements of a basis of $V$ over $k$).
Denote by $\Gal_V$ the Galois group of 
the (separable) closure of $K_V$.
We have canonical maps 
$$
\G^c_K\ra \G^c_V\,\,\, {\rm and}\,\,\,
\G^a_K\ra \G^a_V.
$$ 
If $V$ is 2-dimensional then $K_V$ is isomorphic 
to $k(t)$, (where $1,t$ generate $V$ over $k$).

\begin{lemm}\label{lemm:ms}
Let $V\subset K$ be 2-dimensional.  
Then 
$$
H^2(\Gal_V,\Z/p^n)=0
$$
for all $n\ge 1$. 
\end{lemm}

\begin{proof}
One has
$$
K_2(k(V))=\sum_{\nu} k_{\nu}^*
$$
where the sum is over all codimension 1 points of $k(V)$.
By our assumptions on $k$, the group
$k_{\nu}^*$ is $p$-divisible as well.
Now we apply the theorem of Merkuriev-Suslin (see \cite{MS})
$$
H^2(\Gal_V,\Z/p^n)=K_2(k(V))/p^n=0.
$$
\end{proof}

\begin{rem}
For $K=k(t)$ the irreducible 
divisors are parametrized by $\P^1_k=\P(V)$.
In general, if  $\dim V\ge 3$ then the transcendence
degree of $K_V$ over $k$ is $ > 1$ (and $\le \dim V -1$) and 
the description of $\G^a_V$ is complicated
since there are many more irreducible 
divisors - they cannot be parametrized by an algebraic variety. 
\end{rem}

\subsection{Commuting pairs}  
\label{sect:c-pairs}

If $[\tf_1,\tf_2]=0$ for some lifts to $\G^c$ of elements
$f_1,f_2\in \G^a$ then this commutator vanishes for all 
lifts. In this case we will call the pair $f_1,f_2$
a {\em commuting pair} (c-pair, since the following 
Proposition shows that it is indeed a c-pair in the sense of 
Definition~\ref{defn:c-pair}). 
Natural c-pairs in $\G^a$ arise from
valuations of fields. 

\begin{prop}\label{lemm:cp}
If $f_1,f_2$ are a c-pair then for all 
2-dimensional subspaces 
$V\subset K$ (not necessarily containing $k$)
we have 
$$
\dim \langle f_{1,V},f_{2,V}, 1\rangle \le 2.
$$ 
\end{prop}

\begin{proof}
First assume that $V$ contains $k$.  
We have the following diagram
$$
\begin{array}{ccccccc}
\Gal_K & \ra &  \Gal_K^c & \ra  & \G_K^c & \ra & \G_K^a  \\
\da    &     &    \da    &      &  \da   &     &  \da      \\
\Gal_V & \ra &  \Gal_V^c & \ra  & \G_V^c & \ra & \G_V^a 
\end{array}
$$ 
Lemma~\ref{lemm:ms} implies that
for all surjective continuous homomorphisms
$$
\G^a_V\ra A
$$
onto a finite abelian group $A$
and any cocycle $\al\in H^2(A,\Z/p^n)$ its image in 
$H^2(\Gal_V,\Z/p^n)$ is zero. 
By Lemma~\ref{sect:seq}, it is already zero in 
$H^2(\Gal_V^c,\Z/p^n)$. If $\al$ is nonzero in 
$H^2(Z,\Q/\Z)$ then we can conclude that its image in 
$H^2(\G_V^c,\Z/p^n)$ is zero. 
This means that there exists a
finite group $B$
which is a central extension of $A$ 
and a surjective continuous homomorphism  
$$
\G_V^c\ra B
$$
such that  
$$
H^2(A,\Z/p^n)\ra 0\in H^2(B,Z/p^n).
$$

Assuming that $f_1,f_2$ are nonproportional in 
$\G^a_V$ we construct an $A$ with a nonzero 
cocycle $\al\in H^2(A,\Q/\Z)$ as follows.
Since $f_1,f_2$ are nonproportional 
there exists a sublattice of $k(V)$ of the form 
$\langle x,x+1\rangle$
such that $f_1,f_2$ remain nonproportional on this lattice.
Thus the vectors 
$\hat{f}_1=(f_1(x),f_1(x+1))$ and $\hat{f}_2=(f_2(x),f_2(x+1))$
define a rank 2 lattice $\hat{A}\subset \Z_p\oplus \Z_p$
and we have a surjective homomorphism 
$\G_V^a\ra \hat{A}$.
Then there exists an $n$ such that 
the reduction $A$ 
of $\hat{A}$ modulo $p^n$ is a
subgroup of index $<p^n$ in $\Z/p^n\oplus \Z/p^n$ 
This implies that there 
exists a nonzero cocycle $\al=\al(f_1,f_2)\in H^2(A,\Z/p^n)$
mapping to a nonzero element in $H^2(A,\Q/\Z)$
(by the condition on $\det(\hat{f}_1,\hat{f}_2)$ modulo
$p^n$).
Thus we have surjective maps 
$$
\G_V^c\ra B\ra A
$$ 
where $B$ is a finite group such that $\al$ maps to 
0 in $H^2(B,\Z/p^n)$.

The group 
$B$ contains images of 
$\tilde{f}_1,\tilde{f}_2\in \G_K^c$
and 
$$
\tilde{A}=\langle \tilde{f}_1,\tilde{f}_2\rangle\subset B
$$
surjects onto $A$.  
The cocycle $\al\in H^2(A,\Q/\Z)$ maps to a nonzero element in 
$H^2(\tilde{A},\Q/\Z)$ but to zero in $H^2(B,\Q/\Z)$.  
Contradiction. It follows that the restrictions of 
$f_1,f_2$ to $V$ are proportional.

\

Now we turn to the general case. 
For any 2-dimensional subspace
$V\subset K$ and $x\in V\setminus 0$
consider the 2-dimensional space
$V'$ over $k$ consisting of 
elements of the form $v'=x^{-1}v$ with $v\in V$. 
The space $V'$ contains $k$ and, therefore, 
$f_{1,V'} = \lambda f_{2,V'}$ for some $\la\in k$.
Thus 
$$
f_1(v)+f_1(x^{-1}) = f_1(x^{-1}\,v)=\la f_2(x^{-1}\,v)=
\la(f_2(x^{-1})+f_2(v)),
$$
for all $v \in V$, i.e.,
$\dim\langle f_{1,V},f_{2,V},1\rangle\le 2$.  
\end{proof}

\section{Valuations}

\subsection{Notations}
\label{sect:valuation-notations}

A {\em scale} is a commutative totally ordered group 
(we will use the notations $>$ and $\ge$).  
Let $K$ be a field and 
$$
\nu \,:\ K^*\ra \cI_{\nu}
$$
a surjective homomorphism 
onto a scale $\cI_{\nu}$.  
It is called a nonarchimedean valuation if 
$$
\nu(x+y)\ge \min(\nu(x),\nu(y)),
$$ 
with an equality if $\nu(x)\neq \nu(y)$. 
The group $\cI_{\nu}$ is called the scale of the valuation. 
We consider only nonarchimedean valuations
and call them simply valuations. 

\

A valuation $\nu$ can be extended
to $K$ by $\nu(0)=\infty > \iota$ 
for all $\iota\in \cI_{\nu}$.
It defines a topology on $K$. 
We denote by
$K_{\nu}$ the completion of $K$ with respect to this topology.
The sets  
$$
{\mathcal O}_{\nu,\iota}=\{x\in K\,|\, \nu(x)\ge \iota\}
$$ 
are additive subgroups. 
We call the subring ${\mathcal O}_{\nu}={\mathcal O}_{\nu,\nu(1)}$ 
the {\em valuation ring} of $\nu$.
Denote by 
${\mathfrak m}_\nu= \{x\,|\, \nu(x)> \nu(1)\}$ the 
{\em valuation ideal} (it is a maximal ideal in ${\mathcal O}_\nu$); by  
${\mathcal O}^*_\nu={\mathcal O}_\nu \setminus {\mathfrak m}_\nu$
the set of invertible elements  
and by $\K_\nu={\mathcal O}_\nu/{\mathfrak m}_\nu $ 
the {\em residue field} of ${\mathcal O}_\nu$.  
We have a multiplicative decomposition 
${\mathcal O}^*_\nu/ (1+{\mathfrak m}_\nu)
=\K^*_\nu$.  Here $1+{\mathfrak m}_\nu$ 
is the multiplicative subgroup of ${\mathcal O}^*_\nu$
consisting of elements of the form 
$(1+m)$, $m\in {\mathfrak m}_\nu$.

\subsection{Inertia group}
\label{sect:inertia-group}

Let $K$ and $\nu$ be as above. 
We have a natural embedding 
$\G^a_{K_{\nu}}\ra \G^a_{K}$.
Its image is called the 
{\em abelian valuation group}.

\begin{defn}
Denote by 
$$
\G_{\nu}^a\subset \G^a_K=\LF^{\rm top}(K,\Z_p)
$$ 
the subgroup of those functions which are trivial on 
$(1+{\mathfrak m}_\nu)$. 
This group will be called  
the {\em abelian reduced valuation group}. 
\end{defn}
  
\begin{defn}
Denote by 
$$
{\I}^{a}_{\nu}\subset \G_{\nu}^a\subset \LF^{\rm top}(K,\Z_p)
$$
the subgroup of those functions 
$z^\chi_\nu\in \LF^{\rm top}(K,\Z_p)$ 
such that 
$$
z^\chi_\nu(\kappa)=\chi(\nu(\kappa))
$$
for some homomorphism 
$$
\chi\,:\,  \cI_{\nu}\ra \Z_p
$$
and all $\kappa\in K^*$.
This group will be called  
the {\em  abelian inertia group of $\nu$}. 
The elements $z^\chi_\nu$ are called 
{\em inertia elements} of $\nu$.
\end{defn}

Of course, for all $\kappa,\kappa'\in K^*$ we have
$$
z^\chi_\nu(\kappa\cdot \kappa')=
z^\chi_\nu(\kappa)+z^\chi_\nu(\kappa')
$$
and, since $k$ contains all $p$-power roots,  
we have $\chi(\nu(k^*))=0$.

\begin{rem}
If $\char(\K_\nu)\neq p$ then $\G^a_{\nu}$ 
coincides with the abelian valuation group 
(\cite{ZS}). Otherwise, $\G_{\nu}^a$ is its proper subgroup.
\end{rem}

\subsection{Valuations and flag functions}
\label{sect:val-fg}

\begin{exam}
Let $\nu$ be a valuation on $K=k(X)$ which is trivial 
on $k$. 
Let $z^\chi_\nu$ be an inertia element of $\nu$.
It is a function on $K^*$, invariant under $k^*$. 
We can extend it arbitrarily to $K$, for example by
$z^\chi_\nu(0)=0$.
Then $z^\chi_\nu$ is an abelian flag function on $K$
(considered as a vector space over $k$).
\end{exam}

\begin{exam}
Let $X$ be an algebraic variety defined over $k={\Q}$,
with good reduction $X_p$ at $p$. Let $\nu$ be a divisorial
valuation on the reduction $X_p\otimes \overline{{\mathbf F}}_p$.
This valuation extends to a valuation on $K=\overline{\Q}(X)$, 
with values in $\Q\times \Z$. Any character of $\Z$ which is
trivial on $\Q$ defines an inertia element $z^\chi_{\nu}$, 
which can be considered as a function on $K$.
This is an abelian flag function: for every 
finite dimensional subspace in $K$ the 
corresponding filtration by groups 
consists of modules over 
$p$-integers in $\overline{\Q}$.
\end{exam}

\begin{exam}
Let $K$ be field with a valuation $\nu$,
${\mathcal O}_{\nu}$, ${\mathfrak m}_\nu$, 
${\K}_{\nu}$ as above. Let $\bar{V}$ be 
an $n$-dimensional vector space over ${\K}_{\nu}$
and $\bar{f}$ an AF-function on $\P(\bar{V})$ 
(with respect to ${\K}_{\nu}$).
Define $f$ on ${\mathcal O}_{\nu}^n$, 
extending $f$ trivially over
the cosets ${\mathcal O}_{\nu}/{\mathfrak m}_\nu$. 
Consider the restriction of $f$ to the subset  
$$
V_{\mathcal O}:={\mathcal O}_{\nu}^n\setminus 
({\mathfrak m}_{\nu}{\mathcal O}_{\nu})^n.
$$
Put $V=K^n$ and consider the orbit space $(V\setminus 0)/K^*$. 
Every orbit has a representative in  
$V_{\mathcal O}$. Thus $f$ defines
an AF-function on $\P(V)$ (with respect to $K$).
\end{exam}

\begin{thm}\label{thm:gf-valu} 
Let $f\in \LF(K,\Z_p)\cap \AF(K,\Z_p)$.
Then there exists a valuation $\nu$ on $K$ with scale $\cI_{\nu}$
and a map  $\tilde{f}\,:\, \cI_\nu\ra \Z_p$ such that 
$f(\kappa)=\tilde{f}\circ \nu (\kappa) $ for
all $\kappa\in K$. 
\end{thm}

\begin{proof} 
An AF-function defines a 
filtration $(K^f_{\alpha})_{\al\in \cA}$. 
The logarithmic property of $f$ 
implies that the ordered set $\cA$ is 
an ordered group (a scale). 
The map $\nu\,:\, K\ra \cA$ is a homomorphism. 
Every $K^f_{\alpha}$ is a subgroup under addition. 
Since $f$ is constant on $\ovl{K}^f_{\al}$ it follows that
$\nu$ is a nonarchimedean valuation and that 
$f$ can be factored as claimed.  
\end{proof}

\begin{coro}\label{coro:vvv}
Assume that $f$ satisfies the conditions of 
Theorem~\ref{thm:gf-valu} and let $v$ be the associated valuation.
Consider the groups
$$
\cO_v:=\{ \kappa\in K\,|\, f(\kappa)\ge f(1)\}
$$
$$
\mathfrak m_v := \{ \kappa\in K\,|\, f(\kappa)>f(1)\}.
$$ 
Then $\cO_v/\mathfrak m_v$ is a field and $\cO_v\setminus \mathfrak m_v$
consists of invertible elements in $\cO_v$. 
\end{coro}

\begin{proof}
For any $x\in K^*$
the sets $\overline{K}_f^{\alpha}$ are shifted (bijectively)
under multiplication by $x$. 
In particular, if there is an element $y\in \overline{K}_f^{\alpha}$
such that $xy\in \overline{K}_f^{\alpha}$ then for all
$y'\in  \overline{K}_f^{\alpha}$ the element $xy'$
is also in $\overline{K}_f^{\alpha}$.
The set $\cO_v\setminus \mathfrak m_v$ contains both $x$ and $1\cdot x$
(for any $x\in \cO_v\setminus \mathfrak m_v$).
Thus for all $x\in \cO_v\setminus \mathfrak m_v$ there exists
an inverse in $\cO_v\setminus \mathfrak m_v$.
\end{proof}

\subsection{AF-functions on $\G^a$}

\begin{prop}\label{prop:valua}
Assume that $f_1,f_2\in \G^{a}=\LF^{\rm top}(K,\Z_p)$ 
are linearly independent (over $\Q_p$) 
and that they form a c-pair. Then 
there exists a valuation $\nu$ of $K$ such that
the $\Z_p$-linear span of $f_1,f_2$ contains 
an inertial element $z_\nu^{\chi}\in {\I}_\nu^a$. 
Moreover, for all $\la_1,\la_2\in \Q_p$
the restriction of $\la_1 f_1+\la_2 f_2$
to $1+{\mathfrak m}_v$ is identically 0.
\end{prop}

\begin{proof}
By Proposition~\ref{prop:cp}, 
the $\Z_p$-linear span  of $f_1,f_2$ contains an AF-function.
By Theorem~\ref{thm:gf-valu} and definitions
in Section~\ref{sect:inertia-group}, there exists a 
valuation $\nu$ on $K$ such that this 
AF-function is equal to
$z_{\nu}^{\chi}$ for some inertia element 
$z_{\nu}^{\chi}\in {\I}_{\nu}^a$. 

Let $f$ be any function in the $\Q_p$-linear span of $f_1,f_2$.
Since both $f$ and $z_\nu^{\chi}$ are 
multiplicative, it suffices to consider
them on ${\mathcal O}_v$.
By definition, $z_\nu^{\chi}=0$ on 
$$
{\mathcal O}^*_\nu={\mathcal O}_\nu\setminus {\mathfrak m}_\nu
$$ 
First observe that for $m\in \mathfrak m$ with 
$z_{\nu}^{\chi}(m)\neq 0$ we have $f(1+m)=0$. 
Indeed, consider the sublattice $C=\langle 1,m\rangle$.
Since $z_{\nu}^{\chi}$ is nonconstant on $C$ and 
since it forms a c-pair with $f$, we conclude
that $f$ is proportional to $z_{\nu}^{\chi}$ on this
space. Thus $f(1+m)=z_\nu^{\chi}(1+m)=0$ as claimed. 

Now assume that $z_\nu^{\chi}(m)=0=z_\nu^{\chi}(1)$. 
Then (since $z_\nu^{\chi}$ is an AF-function)
there exists an $m_1\in \mathfrak m$ with 
$z_\nu^{\chi}(m_1)\neq 0$ and $1> m_1>m$. 
Consider the subgroup $\langle m_1,m\rangle$
with generic element $m_1$ and  put $m_2=m_1-m$.  
Then $z_\nu^{\chi}(m')=z_{\nu}^{\chi}(m_2)\neq 0$
and 
$$
f(1+m_1)=f(1-m_2)=0.
$$ 
By Corolllary~\ref{coro:vvv}, we have 
$$
m_3:=\frac{1}{1+m_1-m_2}\in \cO_v\setminus \mathfrak m_v.
$$
Further, by the logarithmic property, 
$$
0 = f(1+m_1)+f(1-m_2)=f(1+m_1-m_2) + f(1-m_1m_2m_3).
$$
Since  
$$
z_\nu^{\chi}(m_1) + z_\nu^{\chi}(m_2)=2 z_\nu^{\chi}(m_1)\neq 0
$$
(as $z_{\nu}^{\chi}$ takes values in $\Z_p$) and 
$z_{\nu}^{\chi}(m_3)=0$  (as $m_3\in \cO_v^*$) we have
that 
$$
f(1-m_1m_2m_3)=0.
$$
This concludes the proof. 
\end{proof}

\begin{coro}
If $f_1,f_2$ satisfy the conditions of \ref{prop:valua} then
there is a valuation $\nu$ 
such that $\langle f_1,f_2\rangle$ lies in
the abelian  reduced valuation group $\G_{\nu}^a$ of $\nu$.
\end{coro}

The subgroup of $\G^a$ generated by $f_1,f_2$ contains a 
cyclic subgroup generated by the
inertial element $z_\nu^{\chi}$ and 
the quotient of $\langle f_1,f_2\rangle$ by
the subgroup generated by AF-elements 
has at most one topological
generator.
An analogous statement is true for liftable abelian
groups of higher rank.

\begin{lemm}\label{lem:18.3}
Let $f_1,\ldots f_n\in \LF(K,\Z_p)$ be 
linearly independent functions. 
Suppose that for every $i,j$ the functions
$f_i, f_j$ form a c-pair.  
Then the group $F$ (topologically) generated by 
$f_1,\ldots f_n$ contains a closed subgroup $F'$
consisting of AF-functions such that 
$F/F'$ is topologically cyclic.
\end{lemm}

\begin{proof}
If all $f_j$ are AF-functions then every $\Z_p$-linear
combination of $f_j$ is an AF-function. 
(Indeed, for every 2-dimensional $V\subset K$
and any pair $f_i,f_j$ the restrictions 
$\rk \,\langle f_{i,V},f_{j,V}\rangle \le 2$. Now apply 
results from \ref{sect:checking})

Assume that $f_1$ is not an AF-function. 
By Proposition~\ref{prop:cp},
there exist $\la_{1j},\la_j\in \Z_p$
such that $\la_{1j}f_1+\la_j f_j$ is an AF-function. 
These functions generate a  $\Z_p$-submodule $F''$ of corank 1, 
such that $F/F''$ is a direct sum of a torsion module and a 
rank one $\Z_p$-module. The torsion elements correspond to AF-functions. 
Denote by $F'$ the module generated by $F''$ and (the preimages
of) these torsion elements. Then $F'$ consists of AF-functions
and $F/F'$ is (topologically)  cyclic.  
\end{proof}

\begin{coro}
The subgroup $F'$ generated by 
AF-elements 
corresponds to the inertia subgroup 
of some valuation $\nu'$.  The lemma implies that every 
liftable noncyclic abelian group lies in some
reduced valuation group and contains 
a group of corank one consisting
of inertial elements.
\end{coro}

\end{document}